\magnification=\magstep 1
\hsize=6.5 true in
\vsize=8.9 true in
\hoffset=0.05 true in

 1

\font\sc=cmcsc10
\font\msbm=msbm10
\font\msbms=msbm7

\tolerance=2000
\hbadness=10000 %never report about underfull horizontal boxes
\hfuzz=2 pt             %after right margin, when needed

\def\frac#1#2{{#1 \over #2}}
\def\frt#1#2{{\textstyle\frac{#1}{#2}}}

\def\qed{\hskip 2em minus 1.5em\vrule height6pt width3pt depth0pt}

\def\eps{\varepsilon}
\def\R{{\hbox{\msbm R}}}
\def\Rs{{\hbox{\msbms R}}}
\def\C{\hbox{\msbm C}}

\def\P{\hbox{\msbm P}\,}

\def\E{\hbox{\msbm E}\,}

\def\cC{{\cal C}}
\def\cE{{\cal E}}
\def\Y{{\cal Y}}
\def\T{{\cal T}}
\def\F{{\cal F}}

\def\A{{\cal A}}

\let\O=\Omega

\let\o=\omega
\def\s{\sigma}
\def\S{\Sigma}
\def\a{\alpha}

\def\l{\lambda}
\def\g{\gamma}

\def\d{\delta}
\def\D{\Delta}
\def\sf{\hbox{$\sigma$-field}}
\let\phi=\varphi

\def\rhol{\rho_{\rm left}}
\def\rhor{\rho_{\rm right}}
\def\U{{\rm U}}
\def\SU{{\rm SU}}
\def\HS{{\rm HS}}
\def\card{{\rm card}\,}
\def\fin{{\rm finite}}

\def\Tr{{\rm tr \,}}

\def\mod{{\rm mod}\,}

\def\equ{$\Longleftrightarrow$}

\def\imp{$\Longrightarrow$}

\def\frt#1#2{{\textstyle\frac{#1}{#2}}}

\def\({\bigl(}
\def\){\bigr)}

\newcount\refno
\refno=0
\newif\iffirst
\firsttrue
\iffirst
\advance\refno by 1
\edef\AKS{\the\refno}
\else
\item{\AKS.} {\sc S.~Albeverio, V.N.~Kolokol'tsov,
O.G.~Smolyanov,} {\it Continuous quantum measurement: local and
global approaches,} {Reviews in Mathematical Physics}
{\bf 9} (1997), 907--920.
\fi

\iffirst
\advance\refno by 1
\edef\AW{\the\refno}
\else
\item{\AW.} {\sc H.~Araki and E.J.~Woods,} {\it Complete Boolean algebras
of type I factors,} Publ.\ RIMS Kyoto Univ., Ser.~A, {\bf 2} (1966),
157--242.
\fi

\iffirst
\advance\refno by 1
\edef\Ar{\the\refno}
\else
\item{\Ar.} {\sc W.~Arveson,} {\it Continuous analogues of Fock
space,} Memoirs AMS {\bf 80}:409 (1989).
\fi

\iffirst
\advance\refno by 1
\edef\Arv{\the\refno}
\else
\item{\Arv.} {\sc W.~Arveson,} {\it $ E_0 $-semigroups in quantum
field theory.} In: ``Quantization, nonlinear partial differential
equations, and operator algebra'', Proc.\
Sympos.\ Pure Math.\ {\bf 59}, AMS 1996, pp.\ 1--26.
\fi

\iffirst
\advance\refno by 1
\edef\Ba{\the\refno}
\else
\item{\Ba.} {\sc P.~Baxendale,} {\it Brownian motions in the
diffeomorphism group I,} {Compositio Mathematica} {\bf 53}
(1984), 19--50. 
\fi

\iffirst
\advance\refno by 1
\edef\Bh{\the\refno}
\else
\item{\Bh.} {\sc B.V.R.~Bhat,} {\it An index theory for quantum
dynamical semigroups,} {Trans.\ Amer.\ Math.\ Soc.} {\bf
348} (1996), 561--583. 
\fi

\iffirst
\advance\refno by 1
\edef\Br{\the\refno}
\else
\item{\Br.} {\sc T.A.~Brun,} {\it Continuous measurements, quantum
trajectories, and decoherent histories.} Report NSF-ITP-97-116,
Inst.\ for Theor.\ Phys., Univ.\ of California, Santa Barbara,
1997. Electronic archive quant-ph/9710021. Submitted to Phys.\ Rev.\ A.
\fi

\iffirst
\advance\refno by 1
\edef\BH{\the\refno}
\else
\item{\BH.} {\sc T.~Byczkowski and A.~Hulanicki,} {\it Gaussian measure of
normal subgroups,} Ann.\ Probab.\ {\bf 11} (1983), 685--691.
\fi

\iffirst
\advance\refno by 1
\edef\BI{\the\refno}
\else
\item{\BI.} {\sc T.~Byczkowski and T.~Inglot,} {\it Gaussian random series
on metric vector spaces,}  Math.\ Zeitschrift {\bf 196}
(1987), 39--50.
\fi

\iffirst
\advance\refno by 1
\edef\CW{\the\refno}
\else
\item{\CW.} {\sc Kai Lai Chung and J.~B.~Walsh,} {\it Meyer's theorem
on predictability,} {Z.~Wahr\-schein\-lich\-keits\-theorie verw.\
Gebiete} {\bf 29} (1974), 253--256.
\fi

\iffirst
\advance\refno by 1
\edef\Da{\the\refno}
\else
\item{\Da.} {\sc E.B.~Davies,} {\sl Quantum theory of open
systems,} Academic Press, London, 1976.
\fi

\iffirst
\advance\refno by 1
\edef\Di{\the\refno}
\else
\item{\Di.} {\sc J.~Dixmier,} {\sl Les $ C^* $-alg\`ebres et leurs
repr\'esentations,} Gauthier-Villars, Paris, 1969.
English translation: {\sl $ C^* $-algebras,} revised edition,
North-Holland 1982.
\fi

% \iffirst
% \advance\refno by 1
% \edef\Du{\the\refno}
% \else
% \item{\Du.} {\sc R.~M.~Dudley,} ``Real analysis and
% probability,'' Wadsworth \& Brooks, Pacific Grove, CA, 1989.
% \fi

\iffirst
\advance\refno by 1
\edef\Fe{\the\refno}
\else
\item{\Fe.} {\sc J.~Feldman,} {\it Decomposable processes and
continuous products of probability spaces,} {J.\ Funct.\
Anal.} {\bf 8} (1971), 1--51.
\fi

\iffirst
\advance\refno by 1
\edef\Ha{\the\refno}
\else
\item{\Ha.} {\sc P.R.~Halmos,} {\sl Lectures on ergodic theory,}
Chelsea, N.Y.\ 1956.
\fi

\iffirst
\advance\refno by 1
\edef\Hu{\the\refno}
\else
\item{\Hu.} {\sc R.L.~Hudson,} {\it Quantum stochastic calculus,
evolutions and flows.} In: ``Quantization, nonlinear partial
differential equations, and operator algebra'', Proc.\
Sympos.\ Pure Math.\ {\bf 59}, AMS 1996, pp.\ 81--91.
\fi

\iffirst
\advance\refno by 1
\edef\Ito{\the\refno}
\else
\item{\Ito.} {\sc K.~It\^o,} {\it Brownian motions in a Lie group,}
Proc.\ Japan Acad., {\bf 26} (1950), 4--10.
\fi

\iffirst
\advance\refno by 1
\edef\ItoS{\the\refno}
\else
\item{\ItoS.} {\sc S.~It\^o,} {\it Brownian motions in a topological group
and in its covering group,} Rend.\ Circ.\ Mat.\ Palermo (2) {\bf 1}
(1952), 40--48.
\fi

\iffirst
\advance\refno by 1
\edef\JK{\the\refno}
\else
\item{\JK.} {\sc N.C.~Jain, G.~Kallianpur,}
{\it Norm convergent expansions for Gaussian processes in Banach spaces,}
Proc.\ Amer.\ Math.\ Soc. {\bf 25} (1970), 890--895.
\fi

\iffirst
\advance\refno by 1
\edef\Ke{\the\refno}
\else
\item{\Ke.} {\sc J.~L.~Kelley,} {\sl General topology,} Berlin, Springer
1955.
\fi

\iffirst
\advance\refno by 1
\edef\Kl{\the\refno}
\else
\item{\Kl.} {\sc V.L.~Klee,} {\it Invariant metrics in groups (solution of
a problem of Banach),} Proc.\ Amer.\ Math.\ Soc.\ {\bf 3} (1952),
484--487.
\fi

\iffirst
\advance\refno by 1
\edef\Ku{\the\refno}
\else
\item{\Ku.} {\sc H.~Kunita,} {\sl Stochastic flows and stochastic
differential equations,} Cambridge Univ.\ Press, 1990.
\fi

\iffirst
\advance\refno by 1
\edef\Li{\the\refno}
\else
\item{\Li.} {\sc G.~Lindblad,}
{\it On the generators of quantum dynamical semigroups,}
{Commun.\ Math.\ Phys.} {\bf 48} (1976), 119--130.
\fi

\iffirst
\advance\refno by 1
\edef\Mal{\the\refno}
\else
\item{\Mal.} {\sc P.~Malliavin,} {\sl Stochastic analysis,}
Springer-Verlag, Berlin, 1997.
\fi

\iffirst
\advance\refno by 1
\edef\Mc{\the\refno}
\else
\item{\Mc.} {\sc H.~P.~McKean,}
{\sl Stochastic integrals,}
Academic Press, New York 1969.
\fi

\iffirst
\advance\refno by 1
\edef\Pa{\the\refno}
\else
\item{\Pa.} {\sc K.R.~Parthasarathy,} {\sl An introduction to
quantum stochastic calculus,} Birk\-h\"auser, Basel, 1992.
\fi

\iffirst
\advance\refno by 1
\edef\Par{\the\refno}
\else
\item{\Par.} {\sc K.R.~Parthasarathy,} {\it Quantum stochastic
calculus,} Proc.\ Intern.\ Congress Math.\ 1994, pp.\
1024--1035, Birkh\"auser, Basel, 1995. 
\fi

\iffirst
\advance\refno by 1
\edef\PZ{\the\refno}
\else
\item{\PZ.} {\sc G.~Da Prato and J.~Zabczyk,} {\sl Stochastic
equations in infinite dimensions,} Cambridge Univ.\ Press 1992.
\fi

\iffirst
\advance\refno by 1
\edef\RS{\the\refno}
\else
\item{\RS.} {\sc M.~Reed and B.~Simon,} {\sl Methods of modern
mathematical physics. I: Functional analysis.}
Revised and enlarged edition, Academic Press, London 1980.
\fi

\iffirst
\advance\refno by 1
\edef\Ro{\the\refno}
\else
\item{\Ro.} {\sc S.~Rolewicz,} {\sl Metric linear spaces,} Warszawa
1972.
\fi

\iffirst
\advance\refno by 1
\edef\Ru{\the\refno}
\else
\item{\Ru.} {\sc T.~de la Rue,} {\it Espaces de Lebesgue,} Lect.\ Notes
Math.\ (Springer) {\bf 1557} (1993), 15--21.
\fi

\iffirst
\advance\refno by 1
\edef\Sa{\the\refno}
\else
\item{\Sa.} {\sc H.\ Sat\^o,} {\it Souslin support and Fourier
expansion of a Gaussian Radon measure,} Lect.\ Notes Math.\ (Springer)
{\bf 860} (1981), 299--313.
\fi

\iffirst
\advance\refno by 1
\edef\Sh{\the\refno}
\else
\item{\Sh.} {\sc A.~Shnirelman,} {\it On the non-uniqueness of weak
solution of the Euler equation.} Preprint IHES/M/96/31.
\fi

\iffirst
\advance\refno by 1
\edef\Sk{\the\refno}
\else
\item{\Sk.} {\sc A.~V.~Skorokhod,} {\sl Asymptotic methods in the
theory of stochastic differential equations,} AMS, Providence,
R.I., 1989 (translated from Russian).
\fi

\iffirst
\advance\refno by 1
\edef\TV{\the\refno}
\else
\item{\TV.} {\sc B.S.~Tsirelson and A.V.~Vershik,} {\it Examples of
nonlinear continuous tensor products of measure spaces and
non-Fock factorizations,} {Reviews in Mathematical Physics}
{\bf 10} (1998), 81--145.
\fi

\iffirst
\advance\refno by 1
\edef\Yo{\the\refno}
\else
\item{\Yo.} {\sc K.~Yosida,} {\it On Brownian motion in a homogeneous
Riemannian space.} Pacific J.\ Math. {\bf 2} (1952), 263--270.
\fi

\hbox{ }

\bigskip\bigskip
\centerline{\sc Unitary Brownian motions are linearizable}
\medskip

\centerline{\sc B.~Tsirelson}
\medskip
\centerline{\it MSRI, Berkeley, and School of Mathematics, Tel Aviv University,}
\centerline{\it Tel Aviv 69978, Israel, e-mail:
 tsirel@math.tau.ac.il}
\bigskip


\medskip

{\narrower
Brownian motions in the infinite-dimensional group of all unitary
operators are studied under strong continuity assumption rather than
norm continuity.  
Every such motion can be described in terms of a countable collection
of independent one-dimensional Brownian motions.
The proof involves continuous tensor products and continuous quantum
measurements.
A by-product: a Brownian motion in a separable F-space (not locally
convex) is a Gaussian process.

}

\bigskip
\centerline{\sc Introduction}
\nobreak\medskip

The most celebrated and useful random process surely is the standard
Brownian motion in $ \R $ (Wiener process). It is Markovian and
Gaussian. Its
increments are independent and stationary. Its continuous sample paths
are bizarre, but many associated probability distributions are smooth,
and connected by wonderful formulas. The multidimensional standard
Brownian motion can produce a lot of random processes by means of
stochastic differential equations. Especially, it can produce its
close relatives, well-known during half a century, --- Brownian
motions in Lie groups and other topological groups.

A Brownian motion in a Lie group $ G $ could be defined
constructively, by means of its generator, an invariant differential
operator of second order on $ G $, or descriptively, as a continuous $
G $-valued random process with stationary independent increments. The
former (constructive) definition stipulates smoothness of the
generator; the latter (descriptive) definition does not. Are they
equivalent? The question was asked by Seizo It\^o, and answered in the
positive by K\^osaku Yosida [\Yo]; see also Kiyosi It\^o [\Ito], and
the book by Henry McKean [\Mc, Sect.~4.7].

Brownian motions in infinite-dimensional groups arise naturally from
stochastic differential equations; see [\Ku, Chap.~4] (``temporally
homogeneous Brownian flows''), [\Mal, Chap.~VIII], [\Sk, Chap.~IV],
[\PZ, Chap.~6]. Does the constructive approach exhaust all
possibilities allowed by the descriptive approach? For the group $ G
= {\rm Diff}^r (M) $ of all diffeomorphisms of a compact smooth
manifold $
M $, the question is answered in the positive by Baxendale [\Ba];
roughly speaking, Brownian motions $ X $ in $ G $ correspond naturally
to Brownian motions $ Y $ in the tangent space $ \D = T_e G $ to $ G $ (at
the unit $ e $ of $ G $) via the stochastic differential equation $
(dX) \cdot X^{-1} = dY $ (in the sense of Stratonovich). For many
other groups, for example, the group of all homeomorphisms of a manifold,
the very idea of $ T_e G $ becomes too vague. Is there any
constructive approach to Brownian motions in such groups?

Let us split the problem in two. First, we waive the relation $ \D =
T_e G $ between $ \D $ and $ G $. Instead, we take the Hilbert space $
\D = l_2 $ once and for all (any other reasonable choice is equivalent,
see 1.9). We try to construct a Brownian motion $ (X,Y) $ in the
direct product $ G \times \D $ (the additive group of $ \D $ is meant),
whose first component $ X $ is the given Brownian motion in $ G $, the
second component $ Y $ is some Brownian motion in $ \D $, and the two
components are perfectly correlated in the following sense: for each $
t \in (0,\infty) $ there is a one-to-one correspondence between sample
paths of $ X $ and $ Y $ on $ [0,t] $; that is, the \sf\ $ \F_t^X $
generated by $ \( X(s) \)_{s\in[0,t]} $ coincides with the \sf\ $
\F_t^Y $ generated by $ \( Y(s) \)_{s\in[0,t]} $. We call such $ (X,Y)
$ a linearization of $ X $.

So, the first part of the problem is linearizability, that is,
existence of a linearization. It
is a well-defined question; either $ X $ is linearizable, or not. If $
X $ is linearizable, we face the second part of the problem: to
find a reasonable interpretation of $ \D $ as a kind of $ T_e G $, and
of the perfect correlation between $ X $ and $ Y $ as a kind of the
stochastic differential
equation $ (dX) \cdot X^{-1} = dY $. (It is not a well-defined
question, in contrast to linearizability.)

The main result of the present work states that every Brownian motion
in the unitary group $ \U(H) $ of the Hilbert space $ H $ (of
countable dimension) is linearizable. The linearizability is evidently
inherited by each group $ G $ that possesses a continuous one-one
homomorphism to $ \U(H) $ (in other words, a faithful unitary
representation); for example, the group of all measure
preserving transformations of $ (0,1) $, or the group of all
diffeomorphisms of a manifold. For the group of all homeomorphisms,
however, the question is still open.

It is meant that the group $ \U(H) $ is equipped with the strong (or,
equivalently, weak) operator topology, rather than the norm topology.
The distinction may be illustrated by the following simple (commutative)
example. Let $ B_k $ be independent standard Brownian motions in $ \R
$. Define a
(random) unitary operator $ X(t) $ by $ X(t) e_k = \exp ( i c_k
B_k(t) ) e_k $, where $ (e_k) $ is an orthogonal basis (thus, $
X(t) $ is diagonal in the basis), and $ (c_k) $ is a sequence of
positive numbers. Then $ X(\cdot) $ is norm continuous for $ c_k = 1 /
\sqrt{\log k} $, but not for $ c_k = 1 $. In the strong operator
topology, however, $ X(\cdot) $ is continuous in any case (even for $
c_k = 2^k $)! Most of well-known results [\Sk], [\PZ] assume norm
continuity, or even stronger conditions.

The crucial idea of splitting the constructive approach in two
separate problems (linearizability, and its interpretation) was
suggested to me by the closely related idea of linearizability in the
theory of continuous tensor products of probability spaces (Feldman
[\Fe], Tsirelson and Vershik [\TV]). I am indebted to Anatoly Vershik
for drawing my attention to the theory in 1994.

The present paper is self-contained; intersections and parallels with
related works are noted, but may be ignored by the reader. Sect.~1
formulates notions and results. Sect.~2 develops a new criterion of
linearizability for continuous tensor products of probability spaces.
Sect.~3 relates unitary Brownian motions to the theory of continuous
quantum measurements (see Davies [\Da]). Sect.~4 establishes local
finiteness of the corresponding quantum stochastic process, which
implies linearizability. It is interesting to observe quantum
probability helping to classical probability. Naturally, driving
forces behind the matter should be more intelligible for readers
acquainted with continuous tensor products of Hilbert spaces and
probability spaces [\AW, \Fe, \Ar, \Arv, \TV], continuous quantum
measurements [\Da, \Li, \Br, \AKS], the relation between the former and the
latter [\Bh], and quantum
stochastic calculus [\Pa, \Par, \Hu].

Brownian motions in a linear topological space are evidently Gaussian,
if the space possesses sufficiently many continuous linear functions.
However, $ L_p $ for $ 0<p<1 $ is an example of a separable 
F-space (complete metric vector space) that possesses no non-zero
continuous linear functions. It
appears that the new criterion of linearizability (Sect.~2) is
fulfilled for all Brownian motions in all linear spaces (in fact, in
all commutative groups), which is shown in Sect.~5. It bridges a gap
between two definitions of Gaussian measures.  (Sect.~5 depends on
Sect.~2, but does not depend on ``quantal'' sections 3,4.)

The following definition, used throughout the paper, is borrowed from
[\ItoS, \Ba].

\medskip

{\sc Definition.} 
A Brownian motion in a topological group $ G $ is a continuous $
G $-valued random process $ X $ such that

(a) $ X(0) = e $ (the unit of $ G $);

(b) increments on the left,
$$  X(t_1), \> X(t_2) (X(t_1))^{-1}, \> \dots, \> X(t_n)
  (X(t_{n-1}))^{-1} \> ,  $$
are independent whenever $ 0 \le t_1 \le \dots \le t_n < \infty $;

(c) the distribution of $ X(t) (X(s))^{-1} $ for $ 0 \le s \le t <
\infty $ depends on $ t-s $ only.

\medskip

For example, any Brownian motion in $ \R $ (the additive group of $
\R $ is meant) is of the form $ X(t) = vt + \s B(t) $ for some $
v \in \R $, $ \s \in [0,\infty) $; here $ B(t) $ is the standard
Brownian motion in $ \R $. Any Brownian motion on the
circe $ \{ z\in\C : |z| = 1 \} $ is of the form $ X(t) = \exp (
ivt + i \s B(t) ) $. The pair $ \( \, \exp ( ivt + i \s B(t) ), \,
B(t) \, \) $
is a Brownian motion in the product (of the circle and the real line),
and forms a linearization of $ \exp ( ivt + i \s B(t) ) $ (for
$ \s \ne 0 $, of course). However, a slightly different terminology is
used in the next section: a Brownian motion in a group is treated as a
representation of a {\it noise} in the group; a linearization of the
motion is treated as a linearization of the noise; and a
linearization of a noise is treated as a faithful representation of
the noise in a linear space.

\bigskip
\centerline{\sc 1. The white noise versus black noises}
\nobreak\medskip

\def\Hlin{H_{\rm lin}}
\def\Flin{\F^{\rm lin}}

Multiplication of operators (or composition of transformations)
gives rise to one-parametric semigroups of operators (or
transformations). Multiplication of measure spaces (or tensor
multiplication of Hilbert spaces) should give rise to
one-parametric semigroups of such spaces. The idea appeared
repeatedly, but still, notions and terminology for spaces
are far less standard than these for operators. ``Continuous
tensor product systems of Hilbert spaces'' are defined by Arveson
[\Ar] in a framework different from that of a theory of
``complete Boolean algebras of type 1 factors'' by Araki and Woods [\AW].
``Factorized Hilbert spaces (and probability spaces) over a
Boolean algebra'' are defined by Tsirelson and Vershik [\TV] in a
framework different from that of Arveson, and of a theory of
``factored probability spaces, indexed by a Borel \sf'' by Feldman
[\Fe]. (The short list of approaches is in no way exhaustive.)
Throughout the paper I restrict myself to the definition given
below, and use the shortest term ``noise'' borrowed from
quantum stochastic calculus.

\medskip

{\sc 1.1 Definition.} A {\it noise\/} consists of a probability
space $ (\O,\F,P) $, a one-parametric group $ (T_t) $ of
measure preserving transformations $ T_t : \O \to \O $ for $ t \in \R $, and a
two-parametric family $ (\F_{s,t}) $ of sub-\sf s $ \F_{s,t} \subset
\F $ for $ -\infty < s \le t < \infty $, such that for all $
r,s,t $

(a) $ T_t $ sends $ \F_{r,s} $ onto $ \F_{r+t,s+t} \, $ ($ r \le
s $);

(b) $ \F_{r,s} $ and $ \F_{s,t} $ are independent ($ r \le s \le
t $);

(c) $ \F_{r,s} $ and $ \F_{s,t} $, taken together, generate $
\F_{r,t} \, $ ($ r \le s \le t $).

\medskip

{\sc 1.2 Note.} Here and henceforth, each probability space is
assumed to be a Lebesgue space (in the sense of Rokhlin, see [\Ru]),
and each \sf\ contains all sets of zero probability.

\medskip

A noise will be called {\it trivial,} if each $ \F_{s,t} $ is
trivial, that is, consists of sets of probability $ 0 $ or $ 1 $
only. Otherwise, each $ \F_{s,t} $ with $ s < t $ is non-atomic.

Remind that a metric (space) is called Polish, if it is complete and
separable. A Polish group is a metrizable topological group $ G $ that
possesses a Polish metric $ \rho $ satisfying the condition
$$  \rho (x_n,y_n) \to 0 \quad \hbox{if and only if both } x_n
  y_n^{-1} \to e \hbox{ and } y_n^{-1} x_n \to e  $$
for all $ x_n, y_n \in G $; here $ e $ is the unit of $ G $. The
condition requires more than $ \rho $ to conform to the topology of $
G $. Every metrizable topological group possesses a left-invariant
metric $ \rhol $ (Birkhoff, Kakutani, see [\Kl]). Existence of a
right-invariant metric follows: $ \rhor (x,y) = \rhol
(x^{-1},y^{-1}) $. In general, no metric is both left-invariant and
right-invariant. Given a left-invariant metric $ \rhol $ and a
right-invariant metric $ \rhor $, we may take $ \rho (x,y) = \rhol
(x,y) + \rhor (x,y) $; we have $ \rhol (x_n,y_n) \to 0 $ \equ\ $
\rhol ( y_n^{-1} x_n, e ) \to 0 $ \equ\ $ y_n^{-1} x_n \to e $;
similarly, $ \rhor (x_n,y_n) \to 0 $ \equ\ $ \rhor (x_n y_n^{-1}, e )
\to 0 $ \equ\ $ x_n y_n^{-1} \to e $, and we get $ \rho (x_n,y_n) \to 0 $
\equ\ ($ y_n^{-1} x_n \to e $ and $ x_n y_n^{-1} \to e $). Thus, $
(x_n) $ is a Cauchy sequence if and only if both $ x_m^{-1} x_n \to e
$ and $ x_n x_m^{-1} \to e $ for $ m,n \to \infty $, which does not
depend on the choice of $ \rhol $, $ \rhor $. So, $ G $ is Polish if
and only if it is complete in the metric $ \rhol + \rhor $ and
separable. See [\Kl], [\Ke, 6-O and 6-Q on pp.~210--213], [\Ba,
Sect.~2].

\medskip

{\sc 1.3 Definition.} Let $ (\O,\F,P) $, $ (T_t) $, $ (\F_{s,t})
$ form a noise, and $ G $ be a Polish group. A {\it
representation} of the noise in the group is a two-parametric
family of $ G $-valued random variables $ X_{s,t}
$ for $ -\infty < s \le t < \infty $, such that for all $ r,s,t
$

(a) $ T_t $ sends $ X_{r,s} $ to $ X_{r+t,s+t} \, $
($ r \le s $),

(b) $ X_{s,t} $ is measurable w.r.t.\ $ \F_{s,t} \, $
($ s \le t $),

(c) $ X_{s,t} X_{r,s} = X_{r,t} \, $
($ s \le t $),

\noindent and for each neighborhoood $ U $ of the unit element of the
group

(d) $ P ( X_{0,t} \in U ) \to 1 $ for $ t \to 0+ $.

\noindent The representation is called {\it continuous,} if for every such
$ U $

(e) $ \frac1t P ( X_{0,t} \notin U ) \to 0 $ for $ t \to 0+ $.

\noindent The representation is called {\it faithful,} if for each $ t \ge
0 $

(f) the \sf\ $ \F_{0,t} $ is generated by the set $ \{ X_{0,s} : s \in
[0,t] \} $ of random variables.

\medskip

{\sc 1.4 Note.} A $ G $-valued random variable is an
equivalence class of measurable maps $ \O \to G $, the
equivalence being equality almost everywhere.

\medskip

It is well-known [\Ba, Th.~3(i)] that the ``continuity
condition'' (e) is equivalent to continuity of $ X_{s,t} $
in $ s $ and $ t $ for almost all $ \o \in \O $ (provided
that functions are appropriately chosen within the
equivalence classes $ X_{s,t} $).

\medskip

{\sc 1.5 Definition.} Let $ G_1, G_2 $ be Polish groups.

(a) $ G_1 $ is {\it Brown subordinate} to $ G_2 $, if every
noise that has a faithful continuous representation in $ G_1 $,
necessarily has a faithful continuous representation in $ G_2
$.

(b) $ G_1 $ and $ G_2 $ are {\it Brown equivalent,} if both $
G_1 $ is Brown subordinate to $ G_2 $, and $ G_2 $ is Brown
subordinate to $ G_1 $.

\medskip

The classical result mentioned in Introduction implies that every
$ n $-dimensional Lie group is Brown equivalent to the $ n
$-dimensional linear space. It is easy to see that an $ m
$-dimensional linear space is Brown subordinate to an $ n
$-dimensional linear space if and only if $ m \le n $.
The proof of the following (main) result is finished at the end of
Sect.~4.

\medskip

{\sc 1.6 Theorem.} The infinite-dimensional unitary group $ \U(H) $ is
Brown equivalent to the Hilbert space (that is, to the
additive group of the Hilbert space with the norm topology).

\medskip

{\sc 1.7 Note.}
Throughout the paper all Hilbert spaces are assumed to be separable
and (unless otherwise stated) of infinite dimension. One Hilbert space,
denoted by $ H $, is complex; it is used as a carrier of the unitary
group $ \U (H) $. Another Hilbert space, denoted by $ \D $, is real; it
is used as a carrier of ``linear'' Brownian motions. The group $ \U
(H) $ is the multiplicative group of all unitary operators on
the Hilbert space $ H $. The strong and the weak operator topologies
coincide on $ \U (H) $; this is the topology $ \U (H) $ is equipped
with. It is metrizable. Here is an example of a right-invariant
metric: $ \rhor (U_1,U_2) = \sum_k 2^{-k} \| U_1 e_k - U_2 e_k \| $;
here $ e_1, e_2, \dots $ are an orthonormal basis of the Hilbert
space. A sequence of unitary operators can converge strongly to a
non-invertible isometric operator, which means that $ \rhor $ is not
complete. However, the metric $ \rho ( U_1, U_2 ) = \rhor (U_1, U_2 )
+ \rhor ( U_1^{-1}, U_2^{-1} ) $ is complete, thus $ \U (H) $ is a
Polish group. (A proof, given in [\Ha, ``Weak topology'' on
pp.~61--64] for the group of all measure preserving transformations of
$ (0,1) $, needs only trivial adaptation to $ \U (H) $.)

\medskip

{\sc 1.8 Theorem.} Every commutative Polish group is Brown
subordinate to the Hilbert space.

\medskip

A separable F-space may be defined as a linear topological space whose
additive group is a Polish group [\Ro, \Kl]. Local convexity is not assumed. The
following facts are evident for locally convex spaces, and well-known
for a number of specific F-spaces (such as $ L_p $ for $ p < 1 $). The
full generality is achieved now by means of the new approach
presented here.

\medskip

{\sc 1.9 Corollary.} All infinite-dimensional separable F-spaces are
Brown equivalent.

\medskip

{\sc 1.10 Corollary} {\it(preliminary formulation).\/} A Brownian motion in a
separable F-space is a Gaussian process.

\medskip

See Sect.~5 for the final formulation of 1.10 (given after discussing
some definitions of Gaussian processes and measures), and for proofs
of 1.8--1.10.

\medskip

{\sc 1.11 Conjecture.} There is a Polish group not Brown subordinate
to the Hilbert space.

\medskip

If a noise has a faithful continuous representation in some
Polish group, then all its representations (in all Polish
groups) are necessarily continuous. That is a consequence of
Meyer's theorem on predictability [\CW] ensuring continuity
of all martingales in the corresponding filtration. We may
avoid using any theory of martingales by means of an equivalent
formulation.
Remind that, given a noise and some $ t > s > 0 $, the
orthogonal projection from $ L_2 (\O,\F_{0,t},P) $ onto $ L_2
(\O,\F_{0,s},P) $ is the conditional expectation, $ f \mapsto
\E ( f | \F_{0,s} ) $. The function $ \E ( f | \F_{0,s} ) $ is
defined up to a negligible set depending on $ s $. We avoid the
trouble of non-countable union of negligible sets by
restricting ourselves to rational numbers $ s $.

\medskip

{\sc 1.12 Definition.} Let $ (\O,\F,P) $, $ (T_t) $, $ (\F_{s,t})
$ form a noise. The noise is called {\it predictable,} if for
any $ t > 0 $ and any $ \F_{0,t} $-measurable bounded function
$ f : \O \to \R $, the function $ \E ( f | \F_{0,s} ) $,
considered for rational $ s \in (0,t) $, is uniformly
continuous in $ s $ for almost all $ \o \in \O $.

\medskip

The word ``predictable'', borrowed from the general theory of
processes and filtrations, does not mean that the future of the
noise can be predicted from its past. It means rather, that
anything is predictable from the infinitesimally near past.
Compare it with the Poisson process; its jumps are utterly
unexpected, they have no precursors.

\medskip

{\sc 1.13 Lemma.} If a noise has a faithful continuous
representation in some Polish group then the noise is
predictable.

\medskip

{\sc 1.14 Lemma.} If a noise is predictable then all its
representations in every Polish group are continuous.

\medskip

Proofs of 1.13 and 1.14 are left to the reader.

\medskip

Predictability is defined in a time-asymmetric way, since $ \E
(f|\F_{0,s}) $ is considered rather than $ \E (f|\F_{s,t}) $.
(You see, $ f \ne \E (f|\F_{0,s}) + \E (f|\F_{s,t}) $ in
general.) Instead of time-reverse predictability, we may ask
about predictability of time-reverse noise, formed by $
(\O',\F',P') = (\O,\F,P) $, $ T'_t = T_{-t} $, and $ \F'_{s,t}
= \F_{-t,-s} $. I do not know, whether predictability of a
noise implies predictability of the time reverse noise, or not.
Also, I do not know, whether every noise is isomorphic to its
time reverse, or not. (A similar question is asked by Arveson
[\Ar, p.~6] about continuous tensor product systems of Hilbert
spaces: ``we do not know if an arbitrary product system must be
antiisomorphic to itself''.) Anyway, if a noise has a faithful
continuous representation in some Polish group then, clearly,
the time-reverse noise has a faithful continuous representation
in the anti-isomorphic group (the same set $ G $, but $ ba $ instead of
$ ab $); both noises are predictable, by Lemma 1.13.

If $ (X_{s,t}) $ is a representation of a predictable noise in
(the additive group of) $ \R $, then $ (X_{0,t})_{t\in[0,\infty)
} $ is a continuous process with stationary independent
increments, that is, a Brownian motion in $ \R $. All such
representations are a linear space, that becomes a (real) Hilbert space
$ \Hlin $ of finite or countable dimension, being equipped with
the norm $ \| (X_{s,t}) \| = \( \E |X_{0,1}|^2 \)^{1/2} $.
Trivial (non-random) representations $ X_{s,t} = v (s-t) $, $ v
\in \R $, form a one-dimensional subspace; its orthogonal
complement $ \Hlin^0 $ consists of centered (that is, zero-mean)
representations; if $ (X_{s,t}) \in \Hlin^0 $ and $ \| (X_{s,t})
\| = 1 $, then $ (X_{0,t})_{t\in[0,\infty)} $ is distributed as
the standard Brownian motion (Wiener process). See [\Ar,
Sect.~5] for the ``dimension of a product system'' parallel to
our $ \dim (\Hlin^0) $. We may choose an orthonormal basis $
(X_{s,t}^k) $, $ 1 \le k < \dim \Hlin $, of $ \Hlin^0 $. If $
\dim \Hlin^0 = d < \infty $, we get a representation $ X_{s,t} =
\( X_{s,t}^1, \dots, X_{s,t}^d \) $ of the noise in $ \R^d $;
otherwise, if $ \dim \Hlin^0 = \infty $, we may construct a
representation $ X_{s,t} = \( c_1 X_{s,t}^1, c_2 X_{s,t}^2,
\dots \) $ in the Hilbert space $ l_2 $, provided that we choose
positive constants $ c_1, c_2, \dots $ decreasing fast enough.
In any case, there is a representation $ (X_{s,t}) $ in the
Hilbert space, containing (in an evident sense) all
representations in $ \R $. For any $ r \le t $ consider the \sf\
$ \Flin_{r,t} $ generated by the set $ \{ X_{r,s} | s \in [r,t]
\} $ of random variables. (The set $ \{ X_{s,t} | s \in [r,t] \}
$ gives the same.) Clearly, $ (\O,\F,P) $, $ (T_t) $, $
(\Flin_{s,t}) $ form a noise, that may be called the {\it
linearizable part} of the given predictable noise $ (\O,\F,P) $,
$ (T_t) $, $ (\F_{s,t}) $.

\medskip

{\sc 1.15 Definition.} A predictable noise is called
{\it linearizable,} if $ \Flin_{s,t} = \F_{s,t} $ for all $ s
\le t $.

\medskip

A predictable noise is linearizable if and only if it has a
faithful (continuous) representation in the Hilbert space.

The following fact results from an example of Tsirelson and
Vershik [\TV, Sect.~5].

\medskip

{\sc 1.16 Theorem.} There is a nontrivial predictable noise with
trivial linear part (that is, having only trivial
representations in $ \R $).

\medskip

The parallel fact for continuous tensor product systems of
Hilbert spaces results from examples of R.T.~Powers, see [\Ar,
Remark 5.4] and [\Arv, p.~12].

So, a linearizable predictable noise may be decomposed into a
finite or countable set of independent copies of the white
noise. The opposite extreme is a predictable noise with trivial
linear part. Such a noise may be called a black noise. Indeed, 
terms ``white noise'' and ``coloured noise'' (these are not
noises as defined here)
imply that a noise manifests itself through linear sensors. For
a black noise, however, the response of any linear sensor is
zero!

What could be a physically reasonable nonlinear sensor able to
sense a black noise? Maybe, a fluid could do it, which is
hinted at by the following words by Shnirelman [\Sh, p.~3]
about a paradoxical motion of an ideal incompressible fluid:
``[\dots] very strong external forces are present, but they are
infinitely-fast oscillating in space, and therefore are
indistinguishable from zero in the sense of distributions.
[\dots] This is the fault of the sensors, not of the forces.''

I do not know, how many nonisomorphic black noises exist, but I
believe they are a continuum (accordingly, the section is
entitled ``the white noise versus black noises'', not ``versus
the black noise''). An invariant, proposed in the next section,
seems to be able to distinguish a continuum of nonisomorphic
black noises. A similar question about continuous
tensor product systems of Hilbert spaces is asked by Arveson
[\Arv, p.~12]: ``It is expected that $ \Sigma $ is uncountable,
but this has not been proved.'' Another question of [\Arv,
p.~12] can be answered (in the positive) by means of Theorem
1.16, namely the question, ``is there a nontrivial $ E_0
$-semigroup $ \a $ with the property that there is a nonzero
unit $ U = \{ U_t : t \ge 0 \} $ and such that every other unit
$ V $ is related to $ U $ by a relation of the form $ V_t =
e^{i\l t} U_t $, $ t \ge 0 $, where $ \l $ is a complex
number?'' In our language, a ``unit'' (for a noise) is a
representation in the multiplicative semigroup of complex
numbers; all such representations are trivial for a black
noise [\TV, Th.~1.7].

\bigskip
\centerline{\sc 2. Spectral type of a noise}
\nobreak\medskip

Given a noise $ (\O,\F,P) $, $ (T_t) $, $ (\F_{s,t}) $, consider
the spaces $ L_2 (\F_{s,t}) = L_2 ( \O, \F_{s,t}, P ) $. If $ r \le
s \le t $ then (under the canonical identification) $ L_2
(\F_{r,s}) \otimes L_2 (\F_{s,t}) = L_2 ( \F_{r,t}) $, which
follows from 1.1(b,c).
At the same time $ L_2 (\F_{r,s}) $ and $ L_2 (\F_{s,t}) $ are
linear subspaces of $ L_2 (\F_{r,t}) $; note that $ L_2 (\F_{r,s})
\cap L_2 (\F_{s,t}) $ is the one-dimensional space of constants,
and $ L_2 (\F_{r,s}) + L_2 (\F_{s,t}) $ is (in general) much
smaller than $ L_2 (\F_{r,t}) $.
Introduce \sf s $ \F_{-\infty,t} $, $
\F_{t,+\infty} $, $ \F_{-\infty,+\infty} $ naturally (say, $
\F_{-\infty,t} $ is generated by all $ \F_{s,t} $ with $ s \in
(-\infty,t] $), then $ \F_{s,t} = \F_{-\infty,t} \cap
\F_{s,+\infty} $ and $ L_2 (\F_{s,t}) = L_2 (\F_{-\infty,t}) \cap
L_2 (\F_{s,+\infty}) $.

Denote by $ E_{s,t} $ the orthogonal projection from $ L_2
(\F_{-\infty,+\infty}) $ onto $ L_2 (\F_{s,t}) $. It is the
conditional expectation: $ E_{s,t} f = \E ( f | \F_{s,t} ) $ for $
-\infty \le s \le t \le +\infty $. Operators $ E_{s,t} $ commute with
each other, and $ E_{s,t} = E_{r,t} E_{s,u} $ whenever $ -\infty
\le r \le s \le t \le u \le +\infty $; in particular, $ E_{s,t} = E_{-\infty,t}
E_{s,+\infty} $.

More generally, given $ -\infty \le s_1 < t_1 < \dots < s_n < t_n
\le +\infty $, we may consider the \sf\ generated by \sf s $
\F_{s_1,t_1}, \dots, \F_{s_n,t_n} $; denote it by $ \F_A $ where $
A = (s_1,t_1) \cup \dots \cup (s_n,t_n) $ is an elementary set,
that is, a finite union of intervals. We identify
elementary sets that coincide up to a finite number of points; say,
$ (r,s) \cup (s,t) $ is identified with $ (r,t) $. More exactly, we
define the elementary Boolean algebra $ \A $ as the factoralgebra
of the Boolean algebra generated by intervals, modulo the
ideal of finite sets. However, it is usual to say ``an elementary
set $ A \in \A $'' instead of ``an equivalence class $ A \in \A
$'', like ``a function $ f \in L_2 $'' instead of ``an equivalence
class $ f \in L_2 $''. For instance, we may say that the complement
of $ (s,t) $ in $ \A $ is $ (-\infty,s) \cup (t,+\infty) $.

So, we have \sf s $ \F_A $ and Hilbert spaces $ L_2 (\F_A) $ for $
A \in \A $, satisfying $ L_2 (\F_{A\cup B}) = L_2 (\F_A) \otimes L_2
(\F_B) $ whenever $ A \cap B = \emptyset $ in $ \A $. (In terms of
[\TV, Def.~1.2] we have a measure factorization over $ \A $.)
Corresponding orthogonal projections $ E_A $ satisfy $ E_A f = \E (
f | \F_A ) $ and $ E_A E_B = E_B E_A = E_{A\cap B} $. However, in
general $ E_A + E_B \ne E_{A\cup B} + E_{A\cap B} $, that is, $
(1-E_A) (1-E_B) \ne 1 - E_{A\cup B} $. The map $ A \mapsto E_A $ is
not a homomorphism of Boolean algebras.

Commuting projections $ E_A $ generate a commutative von Neumann
algebra of operators on the Hilbert space $ L_2
(\F_{-\infty,+\infty}) $; the space decomposes into a direct
integral of Hilbert spaces over the spectrum of the algebra, see
[\Di, Appendix A84]. The spectrum is a (Lebesgue, see 1.2)
measure space $ (Z,\S,\mu) $ (though the
measure is determined up to equivalence), and each $ E_A $ becomes
(the multiplication by) the indicator function $ 1_{e(A)} $ of a set
$ e(A) \in \S $ (or rather, an equivalence class). The relation
$ E_A E_B = E_{A\cap B} $ for operators turns into the relation $
e(A) \cap e(B) = e(A\cap B) $ for sets. So, we have a map,
preserving intersections (but not a
homomorphism) $ A \mapsto e(A) $ from the Boolean algebra $ \A $ to
the Boolean algebra $ \S \, \mod 0 $.

\medskip

{\sc 2.1 Lemma.} For any noise, the \sf\ $ \F_{-\infty,+\infty} $
is generated by the union of \sf s $ \F_{(-\infty,-\eps) \cup
(\eps,+\infty)} $ over all $ \eps > 0 $.

\medskip

{\sc Proof.} In general, an increasing family of \sf s has at most
a countable set of discontinuities (jumps). For the family $
(\F_{-\infty,t})_{t\in\Rs} $ the set of discontinuities must be
shift-invariant due to stationarity, see 1.1(a). Therefore the set
is empty; $ \F_{-\infty,t} $ depends on $ t $ continuously. In
particular, $ \F_{-\infty,0} $ is generated by all $
\F_{-\infty,-\eps} $. Similarly, $ \F_{0,+\infty} $ is generated by
all $ \F_{\eps,+\infty} $. However, $ \F_{-\infty,0} $ and $
\F_{0,+\infty} $ together generate $ \F_{-\infty,+\infty} $.\qed

\medskip

{\sc 2.2 Corollary.} For every $ t \in \R $
$$  \mu \( Z \setminus e ((-\infty,t-\eps) \cup (t+\eps,+\infty)) \) \to 0
  \quad \hbox{for } \eps \to 0 \, .  $$

\medskip

{\sc Proof.} $ E_{(-\infty,-\eps)\cup(\eps,+\infty)} \to 1 $ for $
\eps \to 0 $ strongly on $ L_2 (\F_{(-\infty,+\infty)}) $ by Lemma
2.1, which proves the case $ t = 0 $. The general case follows by
stationarity.\qed

\medskip

Sets $ e(A) $ are determined $ \mod 0 $; to avoid troubles, restrict
ourselves to rational elementary sets $ A $ (I mean that their boundary points
must be rational). For any $ z \in Z $ consider all rational $ A $
such that $ z \in e(A) $. Such $ A $ are a filter (within the
``rational'' subalgebra), since $ e(A) \subset e(B) $ whenever $
A \subset B $, and $ e(A\cap B) = e(A) \cap e(B) $. Though, it may
happen that $ z \in e(\emptyset) $, thus ``filter'' must be understood
here as ``proper or unproper filter''; the unproper filter contains
all $ A $. For each
rational $ t $ there is $ A $ of the filter, bounded away from $ t
$, due to Corollary 2.2. Consider the intersection of all
$ A $ of the filter; boundary points of these $ A $ may be included
or excluded arbitrarily since, being rational, they cannot belong to
the intersection. Denote the intersection by $ C(z) $; it is a
closed set. For every $ t
$, $ \mu \{ z : t \in C(z) \} = 0 $ due to Corollary 2.2. Therefore
$ C(z) $ is of zero Lebesgue measure (for almost all $ z $). Also,
$ C(z) $ is bounded, since $ E_{(-t,t)} \to 1 $ strongly, hence $
\mu \( Z \setminus e((-t,t)) \) \to 0 $ for $ t \to \infty $. So, $
C(z) $ is a nowhere dense compact set.

Sets $ e(A) $ for rational $ A \in \A $ separate points of $ Z $,
therefore $ z \in Z $ is uniquely determined by the corresponding
filter $ \{ A : z \in e(A) \} $. The filter, in its turn, is
uniquely determined by the corresponding intersection $ C(z) $;
namely, a rational $ A \in \A $ belongs to the filter if and only if
$ C(z) $ is contained in the interior of $ A $. So, $ z \mapsto C(z) $ is an
injective map from $ Z $ to the set $ \cC $ of all compact subsets of $ \R
$. The set $ \cC $ becomes a Polish space, being equipped with the
Hausdorff metric. The Borel structure corresponding to the metric
is the same as the Borel structure generated by sets of the
form $ \{ C \in \cC : C \subset A \} $ for all rational $ A \in \A $
(treated as open sets), therefore the map is measurable. We may
identify each $ z $ with $ C(z) $, $ Z $ with $ C(Z) \subset \cC $,
$ \mu $ with a measure on the Polish space $ \cC $ (note that $ \mu (\cC
\setminus Z) = 0 $), and $ \S $ with the \sf\ of $ \mu $-measurable
subsets of $ \cC $, which gives the following result.

\medskip

{\sc 2.3 Theorem.} Let $ (\O,\F,P) $, $ (T_t) $, $ (\F_{s,t}) $
form a noise. Then there are: a probability measure $ \mu $ on the
space $ \cC $ (of all compact subsets of $ \R $, including the empty set),
satisfying the condition
$$  \mu \{ C \in \cC : t \in C \} = 0 \quad \hbox{for all } t \in \R
  \, ,  \eqno \rm(a)  $$
and a direct integral decomposition
$$  L_2 ( \O, \F_{-\infty,+\infty}, P ) = \int_{\cC}^\oplus \hat L (C) \,
  d\mu(C)  \eqno \rm(b)  $$
into a measurable field of Hilbert spaces $ \hat L (C) $
such that for every elementary set (that is, a finite union of
intervals) $ A \subset \R $ and every $ f \in L_2 (\O,\F_{-\infty,+\infty},P) $, decomposed as $ f = \int_{\cC}^{\oplus} \hat f(C)
\, d\mu(C) $, $ \hat f (C) \in \hat L (C) $, the conditional
expectation is decomposed as
$$  \E ( f | \F_A ) = \int_{\cC}^{\oplus} \hat f(C) 1_{e(A)} (C) \, d\mu(C) \,
  ,  \eqno \rm(c) $$
where $ e(A) = \{ C \in \cC : C \subset A \} $.

\medskip

The measure $ \mu $ is determined by the noise up to equivalence.

\medskip

{\sc 2.4 Definition.} The equivalence class of measures $ \mu $ on the
space $ \cC $ (of all compact subsets of $ \R $),
appearing in Theorem 2.3, is called the {\it spectral type} of the
noise. Each such $ \mu $ is called a {\it spectral measure} of the
noise.

\medskip

Given an interval $ (s,t) \subset \R $, we may use the same
construction for decomposing $ L_2 (\F_{s,t}) $ into a direct
integral of spaces $ \hat L_{s,t} (C) $ over the space $ \cC_{s.t} $ of all
compact subsets
of $ (s,t) $ (the subsets are bounded away from $ s $ and $
t $). Given $ r < s < t $, we get two decompositions
of the same Hilbert space,
$$\eqalign{
  & \int^{\oplus} \hat L_{r,t} (C) \, d\mu_{r,t} (C) = L_2 (\F_{r,t}) =
    L_2 (\F_{r,s}) \otimes L_2 (\F_{s,t}) = \cr
  & = \bigg( \int^{\oplus} \hat L_{r,s} (C_1) \, d\mu_{r,s} (C_1) \bigg)
    \otimes \bigg( \int^{\oplus} \hat L_{s,t} (C_2) \, d\mu_{s,t} (C_2)
    \bigg) = \cr
  & = \int^{\oplus} \hat L_{r,s} (C_1) \otimes \hat L_{s,t} (C_2) \, d\mu_{r,s}
    (C_1) d\mu_{s,t} (C_2) \, ,
}$$
which means that for $ \mu_{r,t} $-almost all $ C $
$$  \hat L_{r,t} (C) = \hat L_{r,s} (C\cap(r,s)) \otimes \hat L_{s,t} (C\cap(s,t))
  \leqno (2.5)  $$
(the case $ s \in C $ may be neglected),
and $ \mu_{r,t} = \mu_{r,s} \otimes \mu_{s,t} $. However, 
measures $ \mu_{r,s}, \mu_{s,t}, \mu_{r,t} $ are not
canonical, they are determined up to equivalence;
it is better to write
$$  \mu_{r,t} \sim \mu_{r,s} \otimes \mu_{s,t} \, . \leqno (2.6)  $$
(In terms of [\TV] it is not a measure factorization but a measure
type factorization; the distinction is essential, see the example
at the end of Sect.~1(c) of [\TV]. Another example, closer to
(2.6), is the random set $ C = \{ t : X(t) = 1 \} $, where $ X $ is
the standard Brownian motion in $ \R $; here, $ C \cap (r,s) $ and
$ C \cap (s,t) $ are dependent, but their dependence may be
expressed by a positive density over the product of marginals.)

On the other hand, $ L_2 (\F_{r,s}) $ (as well as $ L_2
(\F_{s,t}) $) is a subspace of $ L_2 (\F_{r,t}) $ (in terms
of the tensor product, $ f $ is identified with $ f \otimes
1 $, where $ 1 $ is the constant function $ 1(\o) = 1 $
treated as a special element of $ L_2 (\F_{s,t}) $), and
the corresponding orthogonal projection $ f \mapsto \E \(
f | \F_{r,s} \) $ transforms $ \int^\oplus_{\cC_{r,t}} \hat
f (C) \, d\mu_{r,t} (C) $ into $ \int^\oplus_{\cC_{r,s}}
\hat f (C) \, d\mu_{r,s} (C) $; here $ \cC_{r,s} $ is
treated as a subset (rather than a factor) of $ \cC_{r,t}
$. Thus, $ \mu_{r,s} \sim \mu_{r,t} \big|_{\cC_{r,s}} $; of
course, the restricted measure $ \mu_{r,t}
\big|_{\cC_{r,s}} $ is defined by $ ( \mu_{r,t}
\big|_{\cC_{r,s}} ) (E) = \mu_{r,t} ( E \cap \cC_{r,s} )
$. Similarly,
$$  \mu_{s,t} \sim \mu \big|_{\cC_{s,t}} \> .  \leqno (2.7)
$$

The spectral measure $ \mu $ on $ \cC $ emerges as follows. Any $ f \in L_2
(\F_{-\infty,+\infty}) $ determines a finite measure $ \mu_f $ on $
\cC $ such that
$$  \mu_f \{ C \in \cC : C \subset A \} = \| \E(f|\F_A) \|^2
  \leqno (2.8)  $$
for all elementary sets $ A \subset \R $. There is $ f $ such that
for every $ g \in L_2 (\F_{-\infty,+\infty}) $, the corresponding
measure $ \mu_g $ is absolutely continuous w.r.t.\ $ \mu_f $ (in
fact, a ``generic'' $ f $ satisfies the condition). For any such $
f $ we may take $ \mu = \mu_f $. Of course, $ \mu_f (E) = \int_E \| \hat f (C) \|^2 \, d\mu(C) $.

In particular, consider the white noise generated by the standard
Brownian motion $ X $ in $ \R $. Any $ f \in L_2 (\F_{-\infty,+\infty})
$ decomposes into multiple It\^o integrals, $ f = \sum_n
\int\!\!\dots\!\!\int \hat f_n (t_1,\dots,t_n) \, dX(t_1) \dots dX(t_n) $.
For the function $ \E ( f | \F_A ) $ the decomposition is the same, but
restricted to $ t_1, \dots, t_n $ belonging to $ A $. The measure $
\mu_f $ is concentrated on finite sets $ C $, and its $ n $-point part
is $ | \hat f_n (t_1,\dots,t_n) |^2 \, dt_1 \dots dt_n $. So, $ \mu $ is
concentrated on finite sets, and its $ n $-point part may be chosen as
the $ n $-dimensional Lebesgue measure.
Spaces $ \hat L (C) $ are one-dimensional, as far as the
Brownian motion $ X $ is one-dimensional; if it is $ d
$-dimensional, then $ \dim \hat L ( \{ t_1,\dots,t_n \} ) = d^n
$.  Note also that the empty set $ C = 
\emptyset \in \cC $ is an atom for $ \mu $, and $ \hat L (\emptyset) $ is the
one-dimensional space of constants, which holds for any noise. We see
that the spectral decomposition of a noise is a generalization of It\^o
decomposition for the white noise.

Consider the set $ \cC_1 \subset \cC $ of all single-element sets $
C $, that is, $ C_1 = \{ \{t\} \big| t \in \R \} $. 
It may happen that $ \mu (\cC_1) > 0 $; in that case we get a nontrivial linear subspace
$$  \hat L_1 = \int_{\cC_1}^{\oplus} \hat L(C) \, d\mu(C) =
\int_\R^\oplus \hat L_1(t) \, dt \subset L_2
  (\F_{-\infty,+\infty}) \, ;  $$
the former integral is the same as in Th.~2.3 but restricted to $
\cC_1 \subset \cC $; it may be transferred to $ \R $ by the one-one
correspondence $ \cC_1 \ni C = \{ t \} \leftrightarrow t \in \R $,
giving the latter integral; the Lebesgue measure ($ dt $) is used,
since the transferred measure is shift-invariant up to equivalence.
Otherwise (when $ \mu(\cC_1)=0 $), $ \hat L_1 $ contains only $ 0 $.
Clearly, $ f \in \hat L_1 $ if and only if $ \mu_f ( \cC \setminus \cC_1
) = 0 $. Each $ f \in \hat L_1 $ gives raise to a family $ (f_{s,t}) $
of $ f_{s,t} \in \hat L_1 \cap L_2 (\F_{s,t}) $ for $ s \le t $ such
that $ \E f_{s,t} = 0 $, and $ f_{r,s} + f_{s,t} = f_{r,t} $ whenever $ r \le s \le t $,
and $ f_{-\infty,+\infty} = f $. (In terms of [\TV] we have an
additive integral, see Def.~1.3 there.) On the other hand, if $ f
\in L_2 (\F_{-\infty,+\infty}) $ belongs to $ L_2 (\F_{-\infty,t}) +
L_2 (\F_{t,+\infty}) $ for every $ t $, and $ \E f = 0 $, then $ f \in \hat L_1 $. (Proof:
$ \mu_f $ is concentrated on $ \{ C \in \cC : \emptyset \ne C \subset (-\infty,t)
\} \cup \{ C \in \cC : \emptyset \ne C \subset (t,+\infty) \} $ for every $ t $;
the intersection over all rational $ t $ gives $ \cC_1 $.) Remind now
the linear part $ (\Flin_{s,t}) $ of a predictable noise, defined in Sect.~1
(before 1.15).

\medskip

{\sc 2.9 Lemma.} For every predictable noise, every $ f \in
\hat L_1 $ is measurable w.r.t.\ $
\Flin_{-\infty,+\infty} $.

\medskip

{\sc Proof.} The space $ \hat L_1 $ is invariant under the
one-parameter unitary group $ (U_t) $ of time shifts,
corresponding to the given group $ (T_t) $ of measure preserving
transformations. Another one-parameter unitary group $ (V_\l) $
acting on $ \hat L_1 $ (but not the whole $ L_2 $) consists of
diagonalizable operators (see [\Di, Appendix A80]) $ V_\l =
\int_\Rs^\oplus e^{i\l t} \, dt $ on $ \hat L_1 = \int_\Rs^\oplus
\hat L_1(t) \, dt $. The two groups satifsy Weil relation $ V_\l U_t =
e^{i\l t} U_t V_\l $. According to the von Neumann uniqueness theorem
(see [\RS, Th.~VIII.14 on p.~275]),
$ \hat L_1 $ decomposes into direct sum
of finite or countable number of irreducible components, ---
subspaces, each carrying an irreducible representation of $ (U_t),
(V_\l) $. Each irreducible representation is unitarily equivalent
to the standard representation in $ L_2 (\R) $, where $ U_t $ acts
as the shift by $ t $, and $ V_\l $ acts as the multiplication by
$ e^{i\l t} $. Comparing the latter with the formula $ V_\l =
\int_\Rs^\oplus e^{i\l t} \, dt $ we conclude that a function $ f
\in L_2 (\R) $ corresponds to $ \int_\Rs^\oplus \hat h(t) f(t) \, dt $
for some vector field $ \hat h(t) \in \hat L_1 (t) $ (not depending on $ f
$). The irreducible component number $ k $ ($ 1 \le k < 1+d $,
where $ d \in \{ 0,1, \dots, \infty \} $ is the number of
components) determines its vector field $ \hat h_k (t) \in \hat L_1 (t) $,
and the set $ \{ \hat h_k(t) | 1 \le k < 1+d \} $ is an orthonormal
basis of $ \hat L_1 (t) $. Comparing the action of $ U_t $ on $ L_2
(\R) $ and $ \hat L_1 $ we conclude that $ U_t \hat h_k(s) = \hat h_k(s+t) $. For
each $ k $ we construct a representation $ (X_{s,t}^k) $ of the
noise in $ \R $ as follows: $ X_{s,t}^k = \int_{(s,t)}^\oplus \hat h_k
(u) \, du $. All $ X_{s,t}^k $ are measurable w.r.t.\ $
\Flin_{-\infty,+\infty} $. Every element of $ \hat L_1 $ is of the form $
\sum_k \int_\R^\oplus \hat h_k (t) f_k (t) \, dt $, therefore it is also
measurable w.r.t.\ $ \Flin_{-\infty,+\infty} $.\qed

\medskip

The same argument can be applied to a non-predictable noise, giving
both Gaussian and Poissonian components of the linear part of the
noise, but we do not need it.

Note that $ \dim \Hlin^0 $, discussed in Sect.~1 (before 1.15), is
equal to $ \dim \hat L_1(t) $, that is, $ \dim \hat L(C) $ for $ C \in \cC_1 $.

\medskip

{\sc 2.10 Corollary.} For every predictable noise, if $ \hat L_1 $
generates the whole \sf\ $ \F_{-\infty,+\infty} $ then the noise is
linearizable.

\medskip

A related result about continuous tensor product systems of Hilbert
spaces is given by Arveson [\Arv, Theorem E in Sect.~6]. His
``decomposable operators'' correspond to a multiplicative counterpart
of $ \hat L_1 $, --- ``multiplicative integrals'' in terms of [\TV], while
elements of $ \hat L_1 $ are ``additive integrals''. The two kinds of
integrals generate the same \sf\ [\TV, Th.~1.7].

\medskip

{\sc 2.11 Corollary.} The following conditions are equivalent for
every predictable noise.

(a) The linear part of the noise is trivial.

(b) $ \Hlin^0 = \{ 0 \} $.

(c) $ \hat L_1 = \{ 0 \} $.

(d) $ \mu (\cC_1) = 0 $.

(e) $ \mu $ is concentrated on sets $ C $ with no isolated points.

\medskip

{\sc Proof.} (a) \equ (b) by definitions; (b) \equ (c) since $ \dim
\Hlin^0 = \dim \hat L_1(t) $; (c) \equ (d) by definition of $ \hat L_1 $; (e) \imp
(d) trivially; and (d) \imp (e) due to (2.7).\qed

\medskip

Turn to the set $ \cC_\fin \subset \cC $ of all finite sets $ C $
(including the empty set). The corresponding subspace

$$  \hat L_\fin = \int^\oplus_{\cC_\fin} \hat L (C) \, d\mu (C)
  \subset L_2 (\F_{-\infty,+\infty})  $$
consists of all $ f \in L_2 (\F_{-\infty,+\infty}) $ such that $
\mu_f ( \cC \setminus \cC_\fin ) = 0 $; the space is non-trivial
if and only if $ \mu (\cC_\fin) > \mu (\{0\}) $.

\medskip

{\sc 2.12 Theorem.} $ \hat L_\fin = L_2 \( \Flin_{-\infty,+\infty}
\) $ for every predictable noise.

\medskip

{\sc Proof.}
$ L_2 \( \Flin_{-\infty,+\infty} \) \subset \hat L_\fin $ due to the
decomposition into multiple Ito integrals, since the linearizable
part of the noise is generated by independent Brownian motions in $
\R $. In order to prove that $ \hat L_\fin \subset L_2 \(
\Flin_{-\infty,+\infty} \) $ note that $ \cC_\fin $ is the union of
sets $ \cC_{t_1,\dots t_n} \subset \cC $ defined for rational $
t_1, \dots, t_n $ such that $ -\infty < t_1 < \dots < t_n < +\infty
$ as follows: $ C \in \cC_{t_1,\dots t_n} $ if and only if each one
of the $ n+1 $ intervals $ (-\infty,t_1), (t_1,t_2),
\dots, (t_{n-1},t_n), (t_n,+\infty) $ contains no more than one point of
$ C $, and no one of the points $ t_1,\dots,t_n $ belongs to $ C $.
Consider the subspace $ \hat L_{t_1,\dots,t_n} = \int_{\cC_{t_1\dots
t_n}}^\oplus \hat L (C) \, d\mu(C) \subset \int_{\cC_\fin}^\oplus \hat L(C) \, d\mu(C) =
\hat L_\fin $. We have $ \cC_{t_1\dots t_n} = \( \{
\emptyset \} \cup (-\infty,t_1) \) \times \( \{ \emptyset \} \cup
(t_1,t_2) \) \times \dots \times \( \{ \emptyset \} \cup (t_n,+\infty)
\) $, where points $ t $ are identified with single-point sets $ \{ t \}
\in \cC $. Thus, $ \hat L_{t_1\dots t_n} =
\( \hat L(\emptyset) \oplus \int_{(-\infty,t_1)}^\oplus \hat L_1(t) \, dt \)
\otimes
\( \hat L(\emptyset) \oplus \int_{(t_1,t_2)}^\oplus \hat L_1(t) \, dt \)
\otimes \dots \otimes
\( \hat L(\emptyset) \oplus \int_{(t_n,+\infty)}^\oplus \hat L_1(t) \, dt \)
$; here $ \hat L(\emptyset) $ is the one-dimensional space of constants.
Combining it with Lemma 2.9 we conclude that all elements of $
\hat L_{t_1,\dots,t_n} $ are measurable w.r.t.\ $ \Flin_{-\infty,+\infty} $.
It remains to note that the union of all $ \hat L_{t_1,\dots,t_n} $
is dense in $ \hat L_\fin $.\qed

\medskip

{\sc 2.13 Corollary.} If $ \hat L_\fin $ generates the whole \sf\ $
\F_{-\infty,+\infty} $, then the noise (assumed to be predictable)
is linearizable.

\medskip

{\sc 2.14 Corollary.} The following conditions are equivalent for
every predictable noise.

(a) A spectral measure is concentrated on finite sets.

(b) The noise is the product of a finite or countable set of
independent copies of the white noise.

\medskip

{\sc 2.15 Note.} Consider the least $ \a $ such that a spectral
measure is concentrated on sets of Hausdorff dimension $ \le \a $. It
is an invariant of a noise. Probably, the invariant takes on a
continuum of values, which could distinguish a continuum of
nonisomorphic black noises.

\bigskip
\centerline{\sc 3. From unitary Brownian motions
  to quantum stochastic processes}
\nobreak\medskip

Brownian motions in a Lie group $ G $ are described by their
generators, right-invariant second-order differential operators on $ G
$. For an $ n $-dimensional $ G $, such an operator depends on $ \frac12
n^2 + \frac32 n $ real parameters. The second-order part of the
operator is given by a symmetric tensor $ (a_{ij}) $ of diffusion
coefficients (at the unit of $ G $), and the first-order part contains
$ n $ more parameters $ (b_i) $, often called the drift vector (at the
unit of $ G $), though they do not form a vector.

Let $ G = \U (d) $ be the group of all unitary $ d \times d $
matrices. The map $ A \mapsto e^{iA} $ from the linear space $ \D $ of
all Hermitian $ d \times d $ matrices to $ \U(d) $ is smooth, and
smoothly invertible in a neighborhood of the origin, which gives a
natural local coordinate system on $ G $; note that $ n = d^2 $.
Denote the inverse map by $ U \mapsto -i \log U $. A Brownian motion $
X $ in $ G $ determines a diffusion process (not a Brownian motion) $
A $ in $ \D $, $ A(t) = -i \log X(t) $, well-defined for small $ t $
(until leaving the neighborhood).

The space $ \D $ is a Euclidean space, with the Hilbert-Schmidt scalar
product $ (A,B)_\HS = \Tr (AB^*) = \Tr (AB) $. Infinitesimal characteristics $
(b_i) $, $ (a_{ij}) $ of $ X $ may be identified with the following
linear and bilinear forms on $ \D $: for all $ B,C \in \D $,
$$  b(B) = \frac{d}{dt} \Big|_{t=0} \E \Tr \( A(t) B \) \> , \quad
  a(B,C) = \frac{d}{dt} \Big|_{t=0} \E \( \Tr (A(t)B) \cdot \Tr (A(t)C)
  \) \> .  $$

There is another interesting parametrization. Matrices $
x_t = \E X(t) $ form a one-parametric semigroup,
therefore $ x_t = e^{tY} $ for some $ Y $; the matrix $
Y $ represents the drift of $ X $ in the linear space of all
matrices (while $ b $ represents the drift of $ X $ in
the manifold of unitary matrices). In order to get 
parameters describing the spread of $
X $, note that each $ U \in \U(d) $ determines an
automorphism $ M \mapsto U M U^* $ of the matrix
algebra $ M_d(\C) $, and the automorphism is a quadratic
function of $ U $. We define
$ \T_t (M) = \E \( X(t) M X^*(t) \) $,
which gives a one-parameter semigroup of linear maps $
\T_t : M_d(\C) \to M_d(\C) $; thus, $ \T_t = e^{tZ} $ for
some linear map $ Z : M_d(\C) \to M_d(\C) $. (Do not confuse
$ \T_t $ with the time shift $ T_t $ introduced by
Def.~1.1.) The following
lemma shows that the pair $ (Y,Z) $ determines $ a $ and $ b $
uniquely. However, the result will not be used, and its
proof is relegated to Appendix.

\medskip

{\sc 3.1 Lemma.} (a) The semigroup $ (\T_t) $ determines uniquely the
law of the Brownian motion $ \( \det X(t) \)^{-1/d} X(t) $ in the
group $ \SU(d) $.

(b) The two semigroups $ (x_t) $, $ (\T_t) $ determine uniquely the
law of the Brownian motion $ X(t) $ in the group $ \U(d) $.

\medskip

Turn to the infinite-dimensional case: $ G = \U(H) $ is the
group of all unitary operators in the Hilbert space (see
1.7). We do not know, how to generalize infinitesimal
characteristics $ a $ and $ b $
for all Brownian motions $ X $ in $ G $ (some cases are
investigated in [\Sk,\PZ]). In contrast, $ Y $ and $ Z $ have a
straightforward generalization. Still, operators 
$$  x_t = \E X(t)  \leqno (3.2)  $$
are well-defined, $ \| x_t \| \le 1 $, and form
a strongly continuous semigroup. The general theory of
operator semigroups ensures that the semigroup has its
generator $ Y $, densely defined, usually unbounded (for
detail, see [\Da]). However, we do not need the generator;
what we need is the very semigroup $ (x_t) $. Another
semigroup $ (\T_t) $ is formed by bounded linear operators
$ \T_t : V \to V $, where $ V $ is the Banach space of all
Hermitian trace-class operators on the Hilbert space $ H $
(see [\Da]); note that $ H $ is complex, while $ V $ is
real. Operators $ \T_t $ are defined as before: 
$$  \T_t (\rho) = \E \( X(t) \rho X^*(t) \)  \leqno (3.3)  $$
for $ \rho \in V $; they form a
strongly continuous semigroup, and $ \| \T_t \| \le 1 $.
In fact, $ (\T_t) $ is a special case of a so-called quantum dynamical
semigroup (see [\Li]), and $ X(\cdot) $ is one of so-called enravelings of $
(\T_t) $ (see [\Br]; physicists are more interested in nonlinear
enravelings). So,
any Brownian motion $ X $ in $ \U(H) $ determines two
semigroups, $ (x_t) $ on $ H $ and $ (\T_t) $ on $ V $. I do
not know, whether $ X $ is uniquely determined by the semigroups,
or not.

Remind the spectral measure $ \mu $ on $ \cC $, and $
\mu_{s,t} $ on $ \cC_{s,t} $. Denote by $ \S $ the \sf\ of
all $ \mu $-measurable subsets of $ \cC $, and by $
\S_{s,t} $ the \sf\ of all $ \mu_{s,t} $-measurable
subsets of $ \cC_{s,t} $.

\medskip

{\sc 3.4 Theorem.} For every Brownian motion in the unitary group $
U(H) $, there is one and only one family $
(\cE_{s,t})_{-\infty<s<t<+\infty} $ satisfying the
following conditions (a)--(f).

(a) For any $ s < t $, $ \cE_{s,t} $ is a map, defined on
the \sf\ $ \S_{s,t} $, taking on values in the algebra
of all bounded linear maps from $ V $ to itself, and $
\cE_{s,t} (E_1) = \cE_{s,t} (E_2) $ whenever $ E_1, E_2 $ differ by
a $ \mu_{s,t} $-negligible set.

(b) The map $ \cE_{s,t} $ is an instrument, as defined in
[\Da, Chap.~4, Def.~1.1]. It means, first, countable additivity:
for any sequence $ E_1, E_2, \dots $ of disjoint sets in $
\S_{s,t} $
$$  \cE_{s,t} ( E_1 \cup E_2 \cup \dots ) \rho = \cE_{s,t}
  (E_1) \rho + \cE_{s,t} (E_2) \rho + \dots  \eqno \rm(b1)  $$
where the sum is norm convergent for each $ \rho \in V $;
second, positivity:
$$  \cE_{s,t} (E) \rho \in V_+ \quad \hbox{for all }
  \rho\in V_+, \, E \in \S_{s,t} \, ;  \eqno \rm(b2)  $$
(here $ V_+ = \{ \rho \in V : \forall h \in H \>\> (\rho h, h) \ge 0 \} $);
and third, trace conservation:
$$  \Tr \( \cE_{s,t}(\cC_{s,t})\rho \) = \Tr (\rho) \quad \hbox{for
  all } \rho \in V \, .  \eqno \rm(b3)  $$
(There is one more condition, complete positivity, see [\Da,
Sect.~9.2]; in fact, it is also satisfied, which, however,
will be neither proved nor used.)

(c) $ \cE_{r,t} ( E_1 \times E_2 ) = \cE_{s,t} (E_2) \cE_{r,s}
(E_1) $ whenever $ r < s < t $, $ E_1 \in \S_{r,s} $, $ E_2 \in
\S_{s,t} $; here $ E_1 \times E_2 $ means $ \{ C \in
\cC_{r,t} : C \cap (r,s) \in E_1 \hbox{ and } C \cap (s,t)
\in E_2 \} $.

(d) $ \cE_{r+t,s+t} ( T_t^{-1} E ) = \cE_{r,s} (E) $
whenever $ r < s $, $ t \in \R $, $ E \in \S_{r,s} $.

(e) $ \cE_{0,t} (\cC_{0,t}) \rho \to \rho $ (in norm) when
$ t \to 0+ $, for every $ \rho \in V $.

(f) $ \cE_{s,t} (\cC_{s,t}) = \T_{t-s} $, and $ \cE_{s,t} \( \{ 
\emptyset \} \) \rho = x_{t-s} \rho
x^*_{t-s} $ whenever $ \rho \in V $ and $ -\infty < s < t <
+\infty $. (Here $ (x_t) $ and $ (\T_t) $ are defined by (3.2), (3.3).)

\medskip

Before starting the proof, remind of some well-known notions and
facts. Let $ \psi \in G \otimes H $, where $ G $ is another
Hilbert space. For any $ h \in H $ define $ (\psi,h)_H \in G $
by the equality $ \( (\psi,h)_H, g \) = \( \psi, g \otimes h \)
$ for all $ g \in G $. For any bounded linear operator $ U : H
\to G \otimes H $ define a bounded linear operator $ \T_U : V
\to V $ as follows. Choose an orthonormal basis $ (g_k) $ in $ G
$, and define $ A_1, A_2, \dots : H \to H $ by $ Uh = \sum_k g_k
\otimes A_k h $ for all $ h \in H $, then $ \T_U (\rho) = \sum_k
A_k \rho A_k^* $. For the one-dimensional operator $ \rho_h \in
V $ defined by $ \rho_h x = (x,h) h $ for all $ x \in H $ we
have $ \( \T_U (\rho_h) h_2, h_1 \) = \sum_k (A_k h, h_1)
\overline{ (A_k h, h_2) } = \( (Uh,h_1)_H, (Uh,h_2)_H \) $ for
all $ h_1, h_2 \in H $, which shows that $ \T_U $ does not
depend on the choice of the basis $ (g_k) $, since linear
combinations of $ \rho_h $ are dense in $ V $.

Let $ F, G $ be Hilbert spaces, and $ U_1 : H \to F \otimes H $,
$ U_2 : H \to G \otimes H $. Define $ U : H \to F \otimes G
\otimes H $ as $ U = ( 1_F \otimes U_2 ) U_1 $. If $ U_1 h =
\sum_k f_k \otimes A_k h $ and $ U_2 h = \sum_l g_l \otimes B_l
h $, then $ U h = \sum_k f_k \otimes U_2 A_k h = \sum_{k,l} f_k
\otimes g_l \otimes B_l A_k h $. It is easy to see that $ \T_U =
\T_{U_2} \T_{U_1} $, that is, $ \T_U (\rho) = \T_{U_2} \(
\T_{U_1} (\rho) \) $ for all $ \rho \in V $.

Let $ U : H \to G \otimes H $ be an isometric operator, and $ G $ be a
direct integral over some measure space, $ G
= \int_\Lambda^\oplus G(\l) \, d\mu(\l) $. For any $ \mu
$-measurable set $ E \subset \Lambda $ consider the projection $
Q_E = \int_\Lambda^\oplus 1_E (\l) \, d\mu(\l) $ on $ G $; that
is, if $ g = \int_\Lambda^\oplus g(\l) \, d\mu(\l) $, then $ Q_E
g = \int_E^\oplus g(\l) \, d\mu(\l) $. Define $ \cE (E) =
\T_{(Q_E \otimes 1_H) U} $. Then $ \cE $ is an instrument.

\medskip

{\sc Proof of Theorem 3.4.} Uniqueness: (c) and (f) determine $
\cE_{s,t}
\( e_{s,t}(A) \) $, where $ e_{s,t} (A) = \{ C \in
\cC_{s,t} : C
\subset A \} $ and $ A \subset \R $ is a finite
union of intervals. Sets $ e_{s,t} (A) $ generate $
\S_{s,t} $ ($ \mod 0 $), and $ e_{s,t} ( A \cap B ) =
e_{s,t} (A) \cap e_{s,t} (B) $. Therefore, $ \cE_{s,t} $ is
determined on the whole $ \S_{s,t} $.

Existence. Introduce an isometric operator $ \tilde X_t : H
\to L_2 (\F_{0,t}, H) $ by $ ( \tilde X_t h ) (\o) = X
(t,\o) h $ for $ h \in H $, $ \o \in \O $. Identifying $
L_2 (\F_{0,t}, H) $ with $ L_2 (\F_{0,t}) \otimes H = \(
\int^\oplus_{\cC_{0,t}} \hat L_{0,t} (C) \, d\mu_{0,t} (C)
\) \otimes H $ we have an isometric operator $ \tilde X_t :
H \to \( \int^\oplus \dots \) \otimes H $.
An instrument $ \cE_{0,t} $ appears, $ \cE_{0,t} (E) =
\T_{(Q_E\otimes1_H) \tilde X_t} $. The same construction for $
X_{s,t}
= X(t) (X(s))^{-1} $ instead of $ X_{0,t} = X(t) $ gives $
\cE_{s,t} $. The construction is invariant under the time
shift group $ (T_t) $. Thus, (a), (b), and (d) hold.

For any $ h_1, h_2 \in H $ we have $ (\tilde X_t h_2, h_1)_H = ( X(t)
h_2, h_1 ) $, since $ ( \tilde X_t h_2, g \otimes h_1 ) = \int_\O \(
X(t,\o) h_2, g(\o) h_1 \) \, dP(\o) = \int_\O \( X(t,\o) h_2, h_1 \)
g(\o) \, dP(\o) $ for all $ g \in L_2 (\F_{0,t}) $.

For any $ \psi \in L_2 (\F_{0,t}) \otimes H $, $ h \in H $, and $ E
\in \S_{0,t} $ we have $ \( (Q_E\otimes1_H) \psi, h \)_H = Q_E
(\psi,h)_H $, since the two expressions are linear in $ \psi $, and
coincide on factorizable vectors; indeed, if $ \psi = f \otimes h_1 $,
$ f \in L_2 (\F_{0,t}) $, $ h_1 \in H $, then both expressions turn
into $ (h_1,h) Q_E f $.

Take $ \psi = \tilde X_t h $, then $ (\psi,h_1)_H = \( X(t)h, h_1 \)
$, thus, $ \( (Q_E\otimes1_H) \tilde X_t h, h_1 \)_H = Q_E \( X(t)h,
h_1 \) $. Applying the formula $ \( \T_U (\rho_h) h_2, h_1 \) = \(
(Uh,h_1)_H, (Uh,h_2)_H \) $ for $ U = (Q_E\otimes1_H) \tilde X_t $, we
get the matrix element
$$  \( \cE_{0,t} (E) (\rho_h) h_2, h_1 \) = \( Q_E ( X(t)h, h_1 ), Q_E
  ( X(t)h, h_2 ) \)  \leqno (3.5)  $$
for all $ h, h_1, h_2 \in H $, $ E \in \S_{0,t} $.

{
For proving (f), consider the two extreme cases, $ E = \{\emptyset\} $ and
$ E = \cC_{0,t} $. For $ E = \cC_{0,t} $ we have $ Q_E =
1_{L_2(\F_{0,t})} $ and $ \( \cE_{0,t} (\cC_{0,t}) (\rho_h) h_2, h_1
\) = \( (X(t)h,h_1), (X(t)h,h_2) \) = \E (X(t)h,h_1) \overline{
(X(t)h,h_2) } $. On the other hand, $ \( X(t)
\rho_h X^*(t) h_2, h_1 \) = ( X^*(t)h_2, h ) ( X(t) h, h_1
) = ( X(t) h, h_1 ) \overline{ ( X(t) h, h_2 ) } $, and we
get $ \( \cE_{0,t} (\cC_{0,t}) (\rho_h) h_2, h_1 \) = \E ( X(t)
\rho_h X^*(t) h_2, h_1 \) = \( \T_t (\rho_h) h_2, h_1 \) $
by (3.3). It means that $
\cE_{0,t} (\cC_{0,t}) = \T_t $. Similarly (or by (e)), $ \cE_{s,t}
(\cC_{s,t}) = \T_{t-s} $, which is the first claim of (f).
For the other case, $ E =
\{ \emptyset \} $, we have $ (Q_{\{ \emptyset \}}f) (\o)
= \E f $ for all $ \o \in \O $, $ f \in L_2 (\F_{0,t}) $,
that is, $ Q_{\{ \emptyset \}} f = ( \E f ) 1_\O $, which is 2.3(c) for
$ A = (-\infty,0) \cup (t,+\infty) $. Thus, $ \( \cE_{0,t} (\{
\emptyset \}) (\rho_h)
h_2, h_1 \) = \( 1_\O \E (X(t)h, h_1), 1_\O \E (X(t)h, h_2) \) = \E
(X(t)h, h_1) \cdot \overline{ \E (X(t)h, h_2) } = ( x_t h, h_1 )
\overline{ ( x_t h, h_2 ) } $.
On the other hand, $ \( x_t \rho_h x^*_t h_2, h_1 \)
= ( x^*_t h_2,h) (x_t h,h_1) = (x_t h,h_1) \overline{ (
x_t h, h_2 ) } $, and we get $ \cE_{0,t} (\{ \emptyset \}) (\rho_h) =
x_t \rho_h x^*_t $ for all $ h \in H $. Therefore $
\cE_{0,t} (\{ \emptyset \}) (\rho) = x_t \rho x^*_t $ for all $ \rho
\in V $. Similarly (or by (d)), $ \cE_{s,t} \( \{ 
\emptyset \} \) \rho = x_{t-s} \rho x^*_{t-s}
$, which completes the proof of (f).
\tolerance=10000\par
}

Item (e) is reduced by (f) to strong continuity of the semigroup $
(\T_t) $, that is, $ \T_t \rho \to \rho $ when $ t \to 0 $. It
suffices to prove it for $ \rho = \rho_h $, $ h \in H $. We have $
\T_t (\rho_h) = \E \( X(t) \rho_h X^*(t) \) = \E \rho_{X(t)h} $;
however, $ \| \rho_{X(t)h} - \rho_h \| \le 2 \| X(t) h - h \| \to 0
$ when $ t \to 0 $ for almost all $ \o $, and never exceeds $ 2 \|
h \| $; therefore $ \| \T_t (\rho_h) - \rho_h \| \le 2 \E \| X(t) h
- h \| \to 0 $, which proves (e).

Before proving (c), note that for any $ r < s < t $, any
linear operator $ U : H \to L_2 (\F_{s,t},H) $ and any
vector-function $ \psi \in L_2 (\F_{r,s},H) $, the
vector-function $ \xi = \( 1_{L_2 (\F_{r,s})} \otimes U \)
\psi $ is given by $ \xi (\o) = \( U(\psi(\o)) \) (\o) $.
Indeed, it suffices to consider a factorizable $ \psi = f \otimes h
$, that is, $ \psi(\o) = f(\o) h $, $ f \in L_2 (\F_{r,s})
$, $ h \in H $. Then $ \xi = f \otimes U h $, that is, $
\xi(\o) = f(\o) \cdot \( (Uh)(\o) \) $, while $ U(\psi(\o))
= U ( f(\o) h ) = f(\o) \cdot (U h) $ and $ \( U(\psi(\o))
\) (\o) = f(\o) \cdot \( (U h) (\o) \) $, too. 

Now assume
that $ r < s < t $, $ E_1 \in \Sigma_{r,s} $, $ E_2 \in
\Sigma_{s,t} $; we have to prove that $ \cE_{r,t}
(E_1\times
E_2) = \cE_{s,t} (E_2) \cE_{r,s} (E_1) $. Similarly to $
\tilde X_{0,t} $, introduce $ \tilde X_{s,t} : H \to L_2
(\F_{s,t},H) = L_2 (\F_{s,t}) \otimes H $ by $ (\tilde
X_{s,t} h) (\o) = X_{s,t} (\o) h $; the same for $ \tilde
X_{r,s} $ and $ \tilde X_{r,t} $. Denote $ F = L_2
(\F_{r,s}) $, $ G = L_2 (\F_{s,t}) $, $ U_1 = \( Q_{E_1}
\otimes 1_H \) \tilde X_{r,s} : H \to F \otimes H $, $ U_2 =
\( Q_{E_2} \otimes 1_H \) \tilde X_{s,t} : H \to G \otimes H
$, and $ U = \( Q_{E_1\times E_2} \otimes 1_H \) \tilde
X_{r,t} : H \to F \otimes G \otimes H $, then (c) takes the
form $ \T_U = \T_{U_2} \T_{U_1} $. It suffices to check that
$ U = ( 1_F \otimes U_2 ) U_1 $. A simpler equality $ \tilde
X_{r,t} = ( 1_F \otimes \tilde X_{s,t} ) \tilde X_{r,s} $
holds, since $ \( ( 1_F \otimes \tilde X_{s,t} )\tilde
X_{r,s} h \) (\o) = \( \tilde X_{s,t} ( (\tilde X_{r,s} h)
(\o) ) \) (\o) = \( \tilde X_{s,t} ( X_{r,s}(\o) h ) \) (\o)
= X_{s,t} (\o) X_{r,s} (\o) h = X_{r,t} (\o) h = (\tilde
X_{r,t} h) (\o) $. Also, $ Q_{E_1\times E_2} = Q_{E_1} \otimes
Q_{E_2} $. Thus, $ U = \( Q_{E_1} \otimes Q_{E_2} \otimes
1_H \) \( 1_F \otimes \tilde X_{s,t} \) \tilde X_{r,s} $,
while $ (1_F\otimes U_2) U_1 = \( 1_F \otimes
(Q_{E_2}\otimes1_H) \tilde X_{s,t} \) \( Q_{E_1}\otimes1_H
\) \tilde X_{r,s} $. It remains to check that $ \( Q_{E_1}
\otimes Q_{E_2} \otimes 1_H \) \( 1_F \otimes \tilde X_{s,t} \)
= \( 1_F \otimes (Q_{E_2}\otimes1_H) \tilde X_{s,t} \)
\( Q_{E_1}\otimes1_H \) $. The two linear operators $ F
\otimes H \to F \otimes G \otimes H $ coincide on the whole
$ F \otimes H $, since they coincide on factorizable
vectors; for any $ f \in F $, $ h \in H $ both operators
transform $ f \otimes h $ into $ Q_{E_1} f \otimes ( Q_{E_2}
\otimes 1_H ) \tilde X_{s,t} h $. So, (c) is verified.\qed

\medskip

The measure $ \mu_f $ on $ \cC_{0,t} $, defined by (2.8) for any $ f
\in L_2 (\F_{0,t}) $, may be written in terms of spectral projections
$ Q_E $ as $ \mu_f (E) = (Q_E f, f) = \| Q_E f \|^2 $. Let $ f $ be a
matrix element of $ X(t) $, that is, $ f = \( X(t) h_2, h_1 \) $ for
arbitrary $ h_1, h_2 \in H $, $ t \in [0,\infty) $. There is a simple
relation between the scalar-valued measure $ \mu_f $ and the
operator-valued measure $ \cE_{0,t} $. The former is a matrix element
of the latter:
$$  \mu_f (E) = \( \cE_{0,t} (E) \rho_{h_2}, \rho_{h_1} \)_\HS \quad
  \hbox{for } f = \( X(t) h_2, h_1 \) \, ;  \leqno (3.6)  $$
here $ (\cdot,\cdot)_\HS $ is the scalar product in the Hilbert space
of all Hilbert-Schmidt operators in $ H $, that is, $ ( \rho_1, \rho_2
)_\HS = \Tr (\rho_1 \rho_2) $ for $ \rho_1, \rho_2 \in V $. (Only
trace-class Hermitian operators are really used here.) Taking into
account that $ (A, \rho_h)_\HS = (Ah,h) $ for every $ A $, we can
deduce (3.6) from (3.5) as follows: $ \( \cE_{0,t} (E) \rho_{h_2},
\rho_{h_1} \)_\HS = \( \cE_{0,t} (E) (\rho_{h_2}) h_1, h_1 \) = \( Q_E
(X(t)h_2,h_1),  Q_E (X(t)h_2,h_1) \) = \| Q_E (X(t)h_2,h_1) \|^2 = \|
Q_E f \|^2 = \mu_f (E) $.

Conditions (a)--(e) of Theorem 3.4 define a notion close to the
notion of a quantum stochastic process as defined in [\Da,
Sect.~5.2]. There are two distinctions. First, Davies stipulates
``the value space $ X $''; assuming that his $ X $ is a single
point, we get rid of it. Second, Davies considers only {\it finite}
sets $ C $, rather than all compact sets. This is the
point! It will be shown in the next section, that our
operator-valued measure $ \cE_{s,t} $ is concentrated on $
\cC_\fin \cap \cC_{s,t} $. Thus, the noise will appear to be linearizable, and $
(\cE_{s,t}) $ will appear to be a quantum stochastic process in the
sense of Davies. However, his ``bounded interaction rate''
condition [\Da, (2.9) in Sect.~5.2] does not hold in our framework;
finiteness of sets $ C $ is ensured by a more subtle mechanism
described in the next section. Note also that the quantum dynamical
semigroup $ (\T_t) $ is strongly continuous ($ \T_t \rho \to \rho $)
but in general not norm continuous; compare it with the following
phrases of Lindblad: ``[\dots] we have to assume that the semigroup is
norm continuous [\dots] a condition which is not fulfilled in many
applications. (We may hope that this restriction can be ultimately
removed using more powerful mathematics.)'' [\Li, p.~120].

\bigskip
\centerline{\sc 4. A compactness argument}
\nobreak\medskip

Let $ X $ be a Brownian motion in the unitary group $ \U(H) $, and
$ (\cE_{s,t}) $ the corresponding family of instruments, given by
Theorem 3.4.
Let $ \rho \in V_1^+ $, that is, $ \rho : H \to H $, $ \rho \ge 0
$, and $ \Tr (\rho) = 1 $. A probability measure $
\mu_{s,t} (\rho,\cdot)
$ on $ ( \cC_{s,t}, \Sigma_{s,t} ) $ arises as follows: $
\mu_{s,t} (\rho,
E) = \Tr \( \cE_{s,t} (E) \rho \) $ for $ E \in \Sigma_{s,t} $.
The measure depends linearly on the parameter $ \rho $. The $ V
$-valued measure $ E
\mapsto \cE_{s,t} (E) \rho $ has its $ V $-valued density $
p_{s,t} (\rho,\cdot) $ w.r.t.\ $ \mu_{s,t} (\rho,\cdot) $; it is a $
\mu_{s,t} (\rho,\cdot) $-measurable function $ \cC_{s,t} \ni C
\mapsto p_{s,t} (\rho,C) \in V_1^+ $
(determined uniquely $ \mu_{s,t} (\rho,\cdot) $-almost everywhere)
such that $ \cE_{s,t} (E) \rho = \int_E p_{s,t} (\rho,C) \,
\mu_{s,t} (\rho,dC) $ for all $ E \in \Sigma_{s,t} $. Its
existence can be checked easily by considering scalar-valued
measures $ E \mapsto \( (\cE_{s,t}(E)\rho) h_2, h_1 \) $ for $
h_1, h_2 \in H $, and their densities (Radon-Nikodym derivatives)
w.r.t.\ $ \mu_{s,t} (\rho,\cdot) $.

{\tolerance=10000
Property 3.4(c), $ \cE_{r,t} (E_1 \times E_2) = \cE_{s,t} (E_2)
\cE_{r,s} (E_1) $, may be reformulated in terms of $ p $ and $
\mu $ as follows. Let $ \rho_0 \in V_1^+ $. Consider $ \rho_1 =
\cE_{r,s} (E_1) \rho_0 = \int_{E_1} p_{r,s} ( \rho_0, C_1 )
\mu_{r,s} ( \rho_0, dC_1 ) $. We have $ \int_{E_1\times E_2}
p_{r,t} ( \rho_0, C ) \mu_{r,t} ( \rho_0, dC ) = \cE_{r,t} (E_1
\times E_2) \rho_0 = \cE_{s,t} (E_2) \rho_1 = \int_{E_1}
\cE_{s,t} (E_2) p_{r,s} ( \rho_0, C_1 ) \mu_{r,s} ( \rho_0, dC_1
) = \int_{E_1} \mu_{r,s} ( \rho_0, dC_1) \int_{E_2} \mu_{s,t} \(
p_{r,s} ( \rho_0, C_1 ), dC_2 \) p_{s,t} \( p_{r,s} ( \rho_0,
C_1 ), C_2 \) $ for all $ E_1 \in \S_{r,s} $, $ E_2 \in \S_{s,t}
$. Taking the trace we get
$$  \mu_{r,t} ( \rho_0, E_1 \times E_2 ) = \int_{E_1} \mu_{r,s}
  ( \rho_0, dC_1 ) \mu_{s,t} \( p_{r,s} ( \rho_0, C_1 ), E_2 \)
  \leqno (4.1)  $$
for $ E_1 \in \S_{r,s} $, $ E_2 \in \S_{s,t} $, and then
$$  p_{r,t} ( \rho_0, C_1 \cup C_2 ) =  p_{s,t} \( p_{r,s} (
  \rho_0, C_1 ), C_2 \) \leqno (4.2)  $$
for almost all $ C_1 \in \cC_{r,s} $, $ C_2 \in \cC_{s,t} $.
The property (4.1) can be generalized to arbitrary $ E \in
\S_{r,t} $; denote $ E(C_1) = \{ C_2 \in \cC_{s,t} : C_1 \cup
C_2 \in E \} $ for $ C_1 \in \cC_{r,s} $, then
$$  \mu_{r,t} ( \rho_0, E ) = \int_{\cC_{r,s}} \mu_{r,s}
  ( \rho_0, dC_1 ) \mu_{s,t} \( p_{r,s} ( \rho_0, C_1 ), E(C_1)
  \) \> .  \leqno (4.3)  $$
It holds by (4.1) for every product set $ E = E_1 \times E_2 $, 
by additivity for finite unions of product sets, and by
continuity for all sets $ E $. Our next step is to make explicit
a Markov property implicit in (4.2), (4.3).

}
Introduce Borel spaces $ \Y_t = \cC_{-\infty,t} \otimes V_1^+ $.
Let $ s < t $. The instrument $ \cE_{s,t} $ gives raise to a
transition probability $ P_{s,t} $, that is, a Borel function on $
\Y_s $ whose values are probability measures on $ \Y_t $. Before
giving a formal definition, consider the idea. We
have $ y_0 \in \Y_s $, that is, $ y_0 = ( C_0, \rho_0 ) $,
$ C_0 \in \cC_{-\infty,s} $ and $ \rho_0 \in V^+_1 $. The
latter determines the corresponding probability distribution $
\mu_{s,t} (\rho_0,\cdot) $ for a random element $ C $ of $
\cC_{s.t} $. Choose $ C $ at random; calculate the
corresponding $ \rho_1 = p_{s,t} (\rho_0,C) \in V_1^+ $.
We get a (random) $ y_1 \in \Y_t $, namely, $ y_1 = ( C_1,
\rho_1 ) $, where $ C_1 = C_0 \cup C $. The distribution
of $ y_1 $ is the measure $ P_{s,t} (y_0) $ on $ \Y_t $
that corresponds to $ y_0 \in \Y_s $. Formally,
$$  P_{s,t} \( (C_0,\rho_0), A \) = \mu_{s,t} \(\rho_0, \{
C \in \cC_{s,t} : (C_0\cup C, p_{s,t} (\rho_0,C) ) \in A
  \} \)  \leqno (4.4)  $$
for all $ C_0 \in \cC_{-\infty,s} $, $ \rho_0 \in V_1^+
$, and Borel sets $ A \subset \Y_t $. The Markov property is
$$  P_{r,t} ( y_0, A ) = \int_{\Y_s} P_{r,s} ( y_0, dy_1 )
  P_{s,t} ( y_1, A )  \leqno (4.5)  $$
for all $ y_0 \in \Y_r $ and all Borel sets $ A \subset \Y_t $.
A proof follows. We have $ y_0 = (C_0,\rho_0) $, $ C_0 \in
\cC_{-\infty,r} $, $ \rho_0 \in V_1^+ $. Introduce $ E = \{ C
\in \cC_{r,t} : \( C_0 \cup C, p_{r,t} (\rho_0,C) \) \in A \} $,
then $ P_{r,t} ( y_0, A ) = \mu_{r,t} ( \rho_0, E ) $ by (4.4),
and $ \mu_{r,t} ( \rho_0, E ) $ is given by (4.3). The
right-hand side of (4.5) is an integral over $ y_1 \in \Y_s $,
but all relevant $ y_1 $ are of the form $ y_1 = \( C_0\cup C_1,
p_{r,s} (\rho_0,C_1) \) $, and the distribution $ P_{r,s}
(y_0,\cdot) $ for $ y_1 $ results from the distribution $
\mu_{r,s} (\rho_0,\cdot) $ for $ C_1 $ by (4.4).
Now (4.5) takes the form
$  \int_{\cC_{r,s}} \mu_{r,s} (\rho_0,dC_1) \mu_{s,t} \(
  p_{r,s} (\rho_0,C_1), E(C_1) \) = \int_{\cC_{r,s}} \mu_{r,s}
  (\rho_0,dC_1) P_{s,t} \( (C_0\cup C_1, p_{r,s} (\rho_0,C_1)),
  A \) $;
it remains to note that $ P_{s,t} \( (C_0\cup C_1, p_{r,s}
(\rho_0,C_1)), A \) = \mu_{s,t} \( p_{r,s} (\rho_0,C_1), E(C_1) \)
$ for $ C_1 \subset \cC_{r,s} $ by (4.4) since, using (4.2),
$ \{ C_2 \in \cC_{s,t} : \( C_0\cup C_1\cup C_2, p_{s,t} ( p_{r,s}
(\rho_0,C_1),C_2) \) \in A \} = \{ C_2 \in \cC_{s,t} : \(
C_0\cup C_1\cup C_2, p_{r,t} (\rho_0,C_1\cup C_2) \) \in A \} =
\{ C_2 \in \cC_{s,t} : C_1\cup C_2 \in E \} = E (C_1) $.

The Markov property (4.5) allows us to introduce a Markov
process $ (Y_t) $. However, having (for now) no information on
its regularity, we restrict ourselves to rational values of $ t
$, thus avoiding the choice of a modification. Strong Markov
property is not claimed, only the simple (fixed-time) Markov
property. Given $ \rho_0 \in V_1^+ $, there is a Markov process
$ (Y_t) $, defined for all rational $ t \in [0,\infty) $, such
that $ Y_t $ takes on values in $ \Y_t = \cC_{-\infty,t} \times
V_1^+ $, that is,
$$\eqalign{
  & Y_t = (C_t,\rho_t), \quad
    Y_0 = ( \emptyset,\rho_0) \hbox{ with probability } 1, \cr
  & \P ( Y_t \in A | Y_s ) = P_{s,t} (Y_s,A) 
    \hbox{ for } s < t, A \subset \Y_t. \cr
}  \leqno (4.6)  $$
The specific form (4.4) of transition probabilities implies
the following. First, $ C_t \cap (-\infty,0) = \emptyset $, and
$ C_t \cap (0,s) $ does not depend on $ t $ as far as $ t \ge s
\ge 0 $. Second, $ \rho_t = p_{s,t} \( \rho_s, C_t \cap (s,t) \)
$ for $ t > s \ge 0 $ (it is meant that $ s,t $ are rational).
Third, $ \P \( C_t \cap (s,t) \in E \big|
C_s, \rho_s \) = \mu_{s,t} ( \rho_s, E ) = \Tr \( \cE_{s,t}(E)
\rho_s \) $ for $ t > s \ge 0 $ and a Borel set $ E \in
\cC_{s,t} $ (note that
the probability does not depend on $ C_s $, and depends linearly
on $ \rho_s $). It follows that $ \E \( \rho_t \big| C_s,\rho_s
\) = \int_{\cC_{s,t}} p_{s,t} (\rho_s,C) \mu_{s,t} (\rho_s,dC) =
\cE_{s,t} (\cC_{s,t}) \rho_s = \T_{t-s} (\rho_s) $ by 3.4(f),
so,
$$  \E \( \rho_t \big| C_s,\rho_s \) = \T_{t-s} (\rho_s)  \leqno
  (4.7)  $$
for rational $ s,t $ such that $ t \ge s \ge 0 $. The other
part of 3.4(f) gives
$$  \P \( C_t \cap (s,t) = \emptyset \big| C_s, \rho_s \) =
  \Tr \( x_{t-s} \rho_s x^*_{t-s} \) \> .  \leqno (4.8)  $$

The space $ V $ has its dual space $ V^* $, identified with the
space of all Hermitian operators $ A : H \to H $; the natural
bilinear form is of course $ A,\rho \mapsto \Tr (A\rho) $. The
semigroup $ (\T_t) $ on $ V $ has its dual semigroup $ \T^*_t $
on $ V^* $; we have $ \Tr \( \rho \T^*_t(A) \) = \Tr \(
\T_t(\rho) A \) = \E \Tr \( X(t) \rho X^*(t) A \) = \E \Tr \(
\rho X^*(t) A X(t) \) $, thus,
$$  \T^*_t (A) = \E \( X^*(t) A X(t) \)  \leqno (4.9)  $$
for $ A \in V^* $ and $ t \ge 0 $. Note that $ V $ is also a
subset of $ V^* $. The following property is peculiar for quantum
dynamical semigroups generated by unitary Brownian motions:
$$  \hbox{If } A \in V \qquad \hbox{then } \T^*_t(A) \in V \quad 
  \hbox{and } \Tr \( \T^*_t(A) \) = \Tr (A) \> ,  \leqno (4.10)
  $$
which follows from (4.9), since $ \Tr \( X^*(t) A X(t) \) = \Tr
\( (X(t))^{-1} A X(t) \) = \Tr (A) $, and $ A \ge 0 $ implies $
X^*(t) A X(t) \ge 0 $.

\medskip

{\sc 4.11 Lemma.} Let $ Q_n \in V^* $ be finite-dimensional
projections, and $ \eps_n \to 0 $ positive numbers. Then the
set $ K = \{ \rho \in V_1^+ : \Tr (Q_n \rho) \ge 1 - \eps_n
\hbox{ for } n=1,2,\dots \} $ is compact.

\medskip

{\sc Proof.} For every $ n $ the set $ \{ Q_n \rho Q_n : \rho
\in V_1^+ \} $ is finite-dimensional and bounded. It suffices
to prove that $ \| \rho - Q_n \rho Q_n \| \le 2 \sqrt{\eps_n}
$ for all $ \rho \in K $. Each $ \rho \in V_1^+ $ may be
represented as $ \rho = \E \rho_h $ for some random vector $
h \in H $, $ \| h \| = 1 $. We have $ \| \rho_h - Q_n \rho_h Q_n
\|^2 = \| \rho_h - \rho_{Q_n h} \|^2 \le 4 \| h - Q_n h \|^2
= 4 \( 1 - \| Q_n h \|^2 \) = 4 \( 1 - \Tr (Q_n \rho_h) \) $,
therefore $ \| \rho - Q_n \rho Q_n \|^2 \le \( \E \| \rho_h -
Q_n \rho_h Q_n \| \)^2 \le \E \| \rho_h - Q_n \rho_h Q_n \|^2
\le 4 \( 1 - \E \Tr (Q_n \rho_h) \) = 4 \( 1 - \Tr (Q_n \rho)
\) \le 4 \eps_n $ for $ \rho \in K $.\qed

\medskip

{\sc 4.12 Lemma.} Let $ S_n \subset V^* $ be compact sets of compact
operators, and $ \eps_n \to 0 $ positive numbers. Then the
set $ K = \{ \rho \in V_1^+ : \max_{A\in S_n} \Tr (A\rho) \ge 1 -
\eps_n \hbox{ for } n=1,2,\dots \} $ is compact.

\medskip

{\sc Proof.} Compactness of an operator $ A \in S_n $ implies
existence of a finite-dimensional projection $ Q_n^{(A)} $
such that $ A \le Q_n^{(A)} + \eps_n \cdot 1_H $. Compactness
of the set $ S_n $ allows to choose a single
finite-dimensional projection $ Q_n $ such that $ A \le Q_n +
\eps_n \cdot 1_H $ for all $ A \in S_n $. Then $ \Tr (A\rho)
\le \Tr (Q_n \rho) + \eps_n $ for all $ A \in S_n $, therefore
$ \Tr (Q_n \rho) \ge \max_{A\in S_n} \Tr (A\rho) - \eps_n \ge
1 - 2 \eps_n $ for all $ \rho \in K $. Lemma 4.11 completes
the proof.\qed

\medskip

{\sc 4.13 Lemma.} For every $ \rho_0 \in V_1^+ $, $ t \in (0,\infty)
$, and $ \eps > 0 $ there is a compact set $ K \subset V_1^+ $ such
that $ \P \( \rho_s \in K \hbox{ for all rational } s \in [0,t] \) \ge 1
- \eps $, where $ \rho_s $ is defined by (4.6).

\medskip

{\sc Proof.} Let $ Q \in V^* $ be a finite-dimensional
projection; introduce a martingale $ M(s) = \E \( \Tr
(Q\rho_t) \big| C_s, \rho_s \) $. Using (4.7), $ M(s) = \Tr \(
Q \T_{t-s} (\rho_s) \) = \Tr ( Q_s \rho_s ) $ where $ Q_s =
\T^*_{t-s} (Q) $. Note that $ Q_s $ is continuous in $ s $
(the same argument as in the proof of 3.4(e)), therefore $ \{
Q_s : s \in [0,t] \} $ is a compact set. Taking into account
that $ M(s) \le 1 $ always, we have $ \P \( \inf_{s\in[0,t]}
M(s) < 1-\eps \) \le \frac1c (1-M(0)) = \frac1c \( 1 - \Tr ( Q
\T_t (\rho_0)) \) $.

Choose finite-dimensional projections $ Q_1, Q_2, \dots $ such
that $ \Tr ( Q_n \T_t (\rho_0) ) \to 1 $. Introduce compact
sets $ S_n = \{ \T^*_{t-s} (Q_n) : s \in [0,t] \} \subset V^*
$; each $ \T^*_{t-s} (Q_n) $ is a compact operator due to
(4.9). Choose $ c_n \to 0 $, $ c_n > 0 $ such that $ \d_n \to
0 $, where $ \d_n = \( 1 - \Tr ( Q_n \T_t (\rho_0)) \) / c_n
$. Martingales $ M_n (s) = \Tr \( \T^*_{t-s} (Q_n) \rho_s \) $
satisfy $ \P \( \inf_{s\in[0,t]} M_n(s) < 1-c_n \) \le \d_n $.
Choosing a subsequence we can get $ \sum \d_n < \infty $; then
$ \P (E_n) \ge 1-\eps_n $, where $ E_n = \{ \o : M_k (s) \ge
1-c_k \hbox{ for all } s \in [0,t] \hbox{ and } k=n,n+1,\dots
\} $ and $ \eps_n = \d_n + \d_{n+1} + \dots \to 0 $. Within $
E_n $ we have $ \Tr \( \T^*_{t-s} (Q_k) \rho_s \) \ge 1-c_k $,
therefore $ \max_{A\in S_k} \Tr (A\rho_s) \ge 1-c_k $ for $
k=n,n+1,\dots $ and all rational $ s \in [0,t] $. Lemma 4.12
ensures that the set $ K_n = \{ \rho \in V_1^+ : \max_{A\in
S_k} \Tr (A\rho) \ge 1-c_k \hbox{ for } k=n,n+1,\dots \} $ is
compact. So, $ \P \( \rho_s \in K_n \hbox{ for all rational }
s \in [0,t] \) \ge \P(E_n) \ge 1-\eps_n \to 1 $.\qed

\medskip

For a given $ \rho_0 \in V_1^+ $ construct compact subsets $
K_1 \subset K_2 \subset \dots $ of $ V_1^+ $ such that $ \P \(
\rho_t \in K_n \hbox{ for all rational } t \in [0,n] \) \ge 1
- 2^{-n} $ for all $ n $. Introduce first exit times $ \tau_n
= \inf \{ t : \rho_t \notin K_n \} $; these are stopping
(Markov) times. (They may be irrational, which is harmless,
since we do not need $ \rho_{\tau_n} $.) We have
$$  \tau_1 \le \tau_2 \le \dots ; \quad \tau_n \to \infty;
  \quad \rho_t \in K_n \hbox{ for } t < \tau_n ; \quad K_1, K_2,
  \dots \hbox{ are compact.}  \leqno (4.14)  $$

Denote $ f_r (\rho) = \Tr \( x_r \rho x^*_r \) = \Tr \(
x^*_r x_r \rho \) $ for $ \rho \in V_1^+ $, $ r \in
[0,\infty) $, then $ \P \( C_{s+r} = C_s \big| C_s, \rho_s \)
= f_r (\rho_s) $. We have $ f_r (\rho) \to 1 $ for $ r \to 0
$, since $ f_r (\rho_h) = \Tr \( \rho_{x_r h} \) = \| x_r h
\|^2 \to \| h \|^2 $ (here $ \rho_h $ is the one-dimensional
operator, $ \rho_h x = (x,h) h $). Functions $ f_r(\cdot) $
are in fact linear functionals on $ V $ of norm $ \le 1 $,
therefore the convergence $ f_r (\cdot) \to 1 $ must be
uniform on compact subsets of $ V_1^+ $. We may choose rational $ r_n
> 0 $ such that $ f_{r_n} (\rho) \ge 1/2 $ for $ \rho \in K_n $; so,
$$  \P \( C_{t+r_n} = C_t \big| C_t, \rho_t \) \ge \frac12 \quad
  \hbox{for } n = 1,2,\dots \hbox{ and all rational } t \in [0,\tau_n)
  \> .  \leqno (4.15)  $$
Consider the number $ \card \( C_{t+r_n} \cap (t,\tau_n) \) $ (maybe,
$ +\infty $) of points in $ C_{t+r_n} $ belonging to the (maybe empty)
interval $ (t,\tau_n) $.

\medskip

{\sc 4.16 Lemma.} For every $ k, n = 1,2,\dots $ and every rational $
t \in [0,\infty) $
$$  \P \( \card \( C_{t+r_n} \cap (t,\tau_n) \) \ge k \) \le 2^{-k} \>
  .  $$

\medskip

{\sc Proof.} Fix $ n $ and $ t $; take an integer $ m $, divide the
interval $ (t, t+r_n) $ into $ m $ intervals $ I_1, \dots, I_m $ of
equal length, and consider the (random) number $ N_m $ of intervals $
I_i $ such that $ I_i \cap C_{t+r_n} \cap (t,\tau_n) \ne \emptyset $.
It suffices to prove that $ \P ( N_m \ge k ) \le 2^{-k} $ for all $ m
$, since $ N_m \to \card \( C_{t+r_n} \cap (t,\tau_n) \) $ for $ m \to
\infty $, and moreover, the convergence is monotone for $ m = 1, 2, 4,
8, \dots $ Further, it suffices to prove that $ \P \( N_m \ge k+1
\big| N_m \ge k \) \le 1/2 $ for $ k=0,1,2,\dots $ However, this fact
follows from (4.15) applied for $ t, t+(1/m)r_n, t+(2/m)r_n, \dots,
t+((m-1)/m)r_n $ by the standard argument with a Markov time, which is
legal, since the Markov time takes on only a finite number of
(rational) values.\qed

\medskip

So, $ C_t \cap [0,\tau_n] $ has a finite intersection with every
interval of length $ r_n $, therefore $ C_t \cap [0,\tau_n] $ is
finite; however, $ \tau_n \to \infty $, thus $ C_t $ is finite for all $ t $ almost sure, that is,
$$  \mu_{0,t} ( \rho_0, \cC_\fin ) = 1 \> .  \leqno (4.17)  $$
We have $ \Tr \( \cE_{0,t} (\cC_\fin) \rho_0 \) = 1 $; $ \Tr \(
\cE_{0,t} (\cC \setminus \cC_\fin) \rho_0 \) = 0 $; by positivity
(3.4b2), $ \cE_{0,t} (\cC \setminus \cC_\fin) \rho_0 = 0 $. However, $
\rho_0 $ is arbitrary; so,
$$  \cE_{0,t} (\cC \setminus \cC_\fin) = 0 \> .  \leqno (4.18)  $$
Combining it with (3.6) we get $ \mu_f (\cC \setminus \cC_\fin) = 0 $
for each $ f $ of the form $ f = \( X(t) h_2, h_1 \) $. It means that
all such $ f $ belong to $ \hat L_\fin $. By Corollary 2.13, the noise
is linearizable.   Thus, the main part of Theorem 1.6 is achieved: the
infinite-dimensional unitary group is Brown subordinate to the Hilbert
space. The converse holds by the two following facts.

\medskip

{\sc 4.19 Note.} Let $ G_1, G_2 $ be Polish groups. If there exists a
continuous one-one homomorphism $ G_1 \to G_2 $, then $ G_1 $ is
Brownian subordinate to $ G_2 $. (The proof is immediate.)

\medskip

{\sc 4.20 Note.} There exists a continuous one-one homomorphism from
(the additive group of) the Hilbert space to the unitary group.

\medskip

{\sc Proof.} Let $ \xi_1, \xi_2, \dots $ be i.i.d.\ $ N(0,1) $ random
variables. Each $ c = (c_1,c_2,\dots) \in l_2 $ determines a random
variable $ \exp \( i \sum c_k \xi_k \) $, and the corresponding
multiplication operator $ U(c) $ on the space of all square integrable
random variables. The operator $ U(c) $ is unitary, and the map $ c
\mapsto U(c) $ is a continuous one-one homomorphism.\qed

\bigskip
\centerline{\sc 5. The commutative case}
\nobreak\medskip

{\sc Proof of Theorem 1.8.}
Let $ G $ be a commutative Polish group, and $ X $ a Brownian motion
in $ G $. Due to Corollary 2.13 it suffices to prove that $ \phi \(
X(t) \) \in \hat L_\fin $ for every bounded Borel function $ \phi : G
\to \R $ and $ t \in [0,\infty) $. That is, we have to prove that $
\mu_f ( \cC_\fin ) = 1 $, where $ f = \phi \( X(t) \) $, $ \| f
\|_{L_2} = 1 $.

We have $ f = \phi ( X_{0,t/2} + X_{t/2,t} ) = \phi ( X_{t/2,t} +
X_{0,t/2} ) $ due to commutatitity of $ G $. It follows that the joint
probability distribution (under $ \mu_f $) of $ C \cap \(0,\frac t2\)
$ and $ C \cap \(\frac t2, t) $ is the same as for $ C \cap \(\frac
t2, t) $ and $ C \cap \(0,\frac t2\) $. That is, $ \mu_f $ is
invariant under the piecewise linear transformation $ \a : (0,t) \to
(0,t) $, $ \a(s) = s + \frac t2 $ for $ s \in \(0,\frac t2\) $, $
\a(s) = s - \frac t2 $ for $ s \in \(\frac t2, t\) $. Likewise, $
\mu_f $ is invariant under the group of all invertible piecewise
linear transformations of $ (0,t) $ having derivative $ =1 $ on each
piece (the group acts naturally on $ \cC_{0,t} $ modulo negligible
sets).

Take an integer $ n $, divide the interval $ (0,t) $ into $ n $
subintervals of length $ t/n $, and consider the probability
distribution (under $ \mu_f $) of the number $ k $ of subintervals
that intersect $ C $; denote the probabilities by $ a_0, a_1, \dots, a_n
$ ($ a_0 + \dots + a_n = 1 $). On the other hand, for any $ \eps \in
(0,1) $ consider $ \d(\eps) = \mu_f \{ C : C \cap [0,\eps t] \ne
\emptyset \} $; we have $ \d(\eps) \to 0 $ for $ \eps \to 0 $ due to 2.3(a). Given
the number $ k $ (of subintervals that intersect $ C $), all the $ n
\choose k $ possibilities are equiprobable due to the invariance
property of $ \mu_f $, and the conditional probability of $ \{ C : C
\cap [0,\frac mn t] = \emptyset \} $ is $ {n-m \choose k} / {n \choose
k} $. So,
$$  \d \bigg( \frac mn \bigg) = 1 - \sum_{k=0}^n a_k \frac{{ n-m
  \choose k }}{{ n \choose k }}  $$
for $ m = 1,\dots,n $. However, $ {n-m \choose k} / {n \choose k} =
\frac{n-m}n \cdot \frac{n-m-1}{n-1} \cdot \dots \cdot
\frac{n-m-k+1}{n-k+1} \le \(\frac{n-m}n\)^k $, therefore $ \d \(\frac
mn\) \ge 1 - \sum a_k \(\frac{n-m}n\)^k = \sum a_k \( 1 -
(\frac{n-m}n)^k \) $, and
$$  a_k + \dots + a_n \le \frac{ \d(\frac mn) }{ 1 - (\frac{n-m}n)^k }
$$
for all $ k $ and $ m $; the left-hand side does not depend on $ m $.

Take an $ \eps \in (0,1) $ and consider $ n = 2,4,8,\dots $ together
with $ m_n $ such that $ m_n / n < \eps $, $ m_n / n \to \eps $ (while
$ k $ does not depend on $ n $). Note that $ a_k + \dots + a_n $
(which should be denoted rigorously by $ a_k^{(n)} + \dots + a_n^{(n)}
$) tends to $ \mu_f \{ C : \card(C) \ge k \} $. We get
$$  \mu_f \{ C : \card(C) \ge k \} \le \frac{ \d(\eps) }{ 1 -
  (1-\eps)^k }  $$
for all $ k=1,2,\dots $ and $ \eps\in(0,1) $; the left-hand side does
not depend on $ \eps $. It remains to choose $ \eps_k \to 0 $ such
that $ (1-\eps_k)^k \le 1/2 $, getting $ \mu_f \{ C : \card(C) \ge k
\} \le 2 \d(\eps_k) \to 0 $, and so, $ \mu_f (\cC_\fin) = 1 $.\qed

\medskip

{\sc Proof of Corollary 1.9.}
A separable F-space is Brown subordinate to the Hilbert space due to
Theorem 1.8; it suffices to prove that the Hilbert space is Brown
subordinate to every infinite-dimensional separable F-space. It
follows immediately from Note 4.19 and the following fact.

\medskip

{\sc 5.1 Note.} Let $ \D $ be the Hilbert space, and $ F $ an
infinite-dimensional F-space. Then there exists a continuous
one-one linear operator $ \D \to F $.

\medskip

{\sc Proof.} Take $ x_1, x_2, \dots \in F $ spanning an
infinite-dimensional subspace. Define $ T : l_2 \to F $ by $ T(\xi_1,
\xi_2, \dots) = \sum \xi_n c_n x_n $, where $ c_n $ tend to $ 0 $ fast
enough, then $ T $ is continuous. Take $ \D = l_2 \ominus T^{-1}(0) $
and restrict $ T $ to $ \D $.\qed

\medskip

The following fact will not be used formally, but is worth to be
mentioned. In particular, it shows that the proof of Corollary 1.10 is
simpler than it seems; path-to-path correspondences used are in fact
point-to-point.

\medskip

{\sc 5.2 Lemma.}
Let $ G $ be a commutative Polish group, $ X $ a Brownian motion in $
G $, $ \D $ the Hilbert space, and $ (X,Y) $ a Brownian motion in
$ G \times \D $ such that $ \F_t^X \subset \F_t^Y $ for all $ t \in
[0,\infty) $. Then for every $ t \in [0,\infty) $, $ X(t) $ is
measurable w.r.t.\ the \sf\ generated by $ Y(t) $.

\medskip

{\sc Proof.}
It suffices to prove that the random variable $ f = \phi(X(t)) $ is
measurable w.r.t.\ the
\sf\ generated by $ Y(t) $ for every bounded Borel function $ \phi : G
\to \R $. We may assume that $ \D $ is a space of sequences, and $ Y(t)
= \( B_1(t), B_2(t), \dots \) $, where $ B_1, B_2, \dots $ are
independent standard Brownian motions (in $ \R $). Like every element
of $ L_2 (\F_t^Y) $, $ f $ can be decomposed into multiple It\^o
integrals,
$$  f = \sum_{n=0}^\infty \int\!\!\dots\!\!\int \sum_{k_1,\dots,k_n}
  \hat f_{k_1,\dots,k_n} (t_1,\dots,t_n) \, dB_{k_1}(t_1) \dots
  dB_{k_n}(t_n) \> .  $$
Remind the invariance argument used in the proof of Theorem 1.8:
piecewise linear transformations of $ (0,t) $ (with derivative $ =1 $
on each piece) act on paths of $ (X,Y) $ by measure preserving
transformations, leaving invariant $ X(t) $ and $ f $. Therefore, $
\hat f_{k_1,\dots,k_n} (t_1,\dots,t_n) $ remains unchanged, when the
transformation acts on each of $ t_1,\dots,t_n $. It means that $ \hat
f_{k_1,\dots,k_n} (t_1,\dots,t_n) $ does not depend on $ t_1,\dots,t_n
$, it is a constant! Thus, the stochastic integral is a polynomial of
$ Y(t) $ only.\qed

\medskip

Denote by $ \g $ the standard Gaussian measure on the space $
\R^\infty $ of all sequences of reals; that is, $ \g $ is the joint
distribution of a sequence of i.i.d.\ $ N(0,1) $ random variables.
Given a Polish group $ G $, we introduce the set $ L_0 (\g,G) $ of all
equivalence classes of $ \g $-measurable functions $ X : \R^\infty \to
G $, the equivalence being the equality $ \g $-almost everywhere.

\medskip

{\sc 5.3 Lemma.}
Let $ G $ be a commutative Polish group.
For any $ X \in L_0 (\g,G) $ the following three properties are
equivalent.

(a) There exists $ Y \in L_0 (\g,G) $ such that $ X(u) + X(v) = Y \(
2^{-1/2} (u+v) \) $ for $ \g\otimes\g $-almost all pairs $ (u,v) \in
\R^\infty \times \R^\infty $.

(b) For each $ a \in l_2 $ there exists $ g \in G $ such that $ X(u+a)
= X(u) + g $ for $ \g $-almost all $ u \in \R^\infty $.

(c) There exist a subgroup $ E \subset \R^\infty $ of $ \g $-full
measure, a homomorphism $ X_1 : E \to G $, and an element $ g \in G $
such that $ X(u) = X_1 (u) + g $ for $ \g $-almost all $ u $.

\medskip

{\sc Proof.}
(c) \imp (a): $ X(u) + X(v) = X_1 (u) + g + X_1 (v) + g = X_1 (u+v) +
2g $; we let $ Y(w) = X_1 (2^{1/2} w) + 2g $, taking into account that $
2^{1/2} w \in E $ for $ \g $-almost all $ w $.

(a) \imp (b): $ X(u+a) + X(v-a) = Y \( 2^{-1/2} (u+v) \) = X(u) + X(v)
$, therefore $ X(u+a) - X(u) = X(v) - X(v-a) $ for $ \g \otimes \g
$-almost all $ (u,v) $, which means that both functions are constant $
\g $-almost everywhere.

(b) \imp (c).
The function $ X_0 : l_2 \to G $ defined by $ X(u+a) = X(u) + X_0(a)
$ is a homomorphism. If $ a_1, a_2, \dots \in l_2 $, $ a_n \to 0 $
weakly in $ l_2 $, then $ X_0 (a_n) \to 0_G $ for $ n\to\infty $,
since $ X (\cdot+a_n) \to X(\cdot) $ in probability. Introduce
projections $ P_n : \R^\infty \to \R^n \subset l_2 \subset \R^\infty
$, $ P_n (\xi_1,\xi_2,\dots) = (\xi_1,\dots,\xi_n) =
(\xi_1,\dots,\xi_n, 0,0,\dots) $, and functions $ X_n = X - X_0
\circ P_n $, that is, $ X_n (\xi_1,\xi_2,\dots) = X
(\xi_1,\xi_2,\dots) - X_0 (\xi_1,\dots,\xi_n, 0,0,\dots) $, then $ X_n
(u+a) = X_n (u) $ for all $ a \in \R^n $, which means that $ X_n
(\xi_1,\xi_2,\dots) $ depends only on $ \xi_{n+1}, \xi_{n+2}, \dots $

Functions of the form $ \phi_n \circ P_n $ (for all $ n $ and all
measurable $ \phi_n : \R^n \to G $) are dense in $ L_0 (\g,G) $
(equipped with the convergence in probability). Choose $ \phi_n $ such
that $ \phi_n \circ P_n \to X $ in probability. We have $ \phi_n \circ
P_n - X \to 0 $ in probability. Choose $ \eps_n \to 0 $
such that $ \g \{ u : \rho \( 0, \phi_n(P_n u) - X(u) \) \le \eps_n \}
\ge 1 - \eps_n $. The probability of the event $ \rho \( 0, \phi_n(P_n
u) - X(u) \) \le \eps_n $ is the expectation of the conditional
probability of the same event, given the first $ n $ coordinates $
\xi_1,\dots,\xi_n $ of $ u = (\xi_1,\xi_2,\dots) $. Given $ n $, we
may choose $ \xi_1,\dots,\xi_n $ such that the conditional probability
at $ (\xi_1,\dots,\xi_n) $ is $ \ge 1 - \eps_n $. However, $ \phi_n
(P_n u) - X(u) = \phi_n(P_n u) - X_0(P_n u) - X_n (u) $. Denote $ g_n
= \phi_n (\xi_1,\dots,\xi_n) - X_0 ( \xi_1,\dots,\xi_n, 0,0,\dots )
$, then $ \g \{ u : \rho \( 0, g_n - X_n(u) \) \le \eps_n \} \ge 1 -
\eps_n $; the conditioning is omitted, since $ X_n (\xi_1,\xi_2,\dots) $
depends only on $ \xi_{n+1}, \xi_{n+2}, \dots $
So, there exist $ g_1,g_2,\dots \in G $ such that $ X_n - g_n \to 0 $
in probability, which means that $ X - X_0 \circ P_n - g_n \to 0 $,
that is, $ X_0 \circ P_n + g_n \to X $ in probability, when $
n\to\infty $. The following trick shows that $ g_n \to g $ for some $
g $.

Consider three measure preserving maps $ \a,\a_1,\a_2 : \( \R^\infty
\times \R^\infty, \g\otimes\g \) \to ( \R^\infty, \g ) $ defined by $
\a_1 (u,v) = \frac12 u + \frac{\sqrt3}2 v $, $ \a_2 (u,v) = \frac12 u -
\frac{\sqrt3}2 v $, $ \a (u,v) = u $, then $ \a_1(u,v) + \a_2(u,v) =
\a(u,v) $. For any $ a,b \in l_2 $ we have $ (X\circ\a_1 + X\circ\a_2 -
X\circ\a) (u+a,v+b) - (X\circ\a_1 + X\circ\a_2 - X\circ\a) (u,v) =
X_0(\a_1(a,b)) + X_0(\a_2(a,b)) - X_0(\a(a,b)) = X_0 \( \a_1(a,b) +
\a_2(a,b) -
\a(a,b) \) = 0 $, which means that the function $ X\circ\a_1 +
X\circ\a_2
- X\circ\a $ is constant; $ X(\a_1(u,v)) + X(\a_2(u,v)) - X(\a(u,v)) = g $
for $ \g\otimes\g $-almost all $ (u,v) $. However, $ X \circ \a_1 =
\lim_n (X_0\circ P_n\circ\a_1 + g_n) $, and the same for $ \a_2,\a $. We
have $ g = \lim_n ( X_0\circ P_n\circ\a_1 + g_n + X_0\circ
P_n\circ\a_2 +
g_n - X_0\circ P_n\circ\a - g_n ) = \lim_n ( X_0\circ
P_n\circ(\a_1+\a_2-\a) + g_n ) = \lim_n g_n $.

So, $ X_0\circ P_n + g \to X $ in probability, when $ n\to\infty $. We
choose $ n_k \to \infty $ such that $ X_0 \circ P_{n_k} + g \to X $
almost sure. The set of all $ u $ such that $ \lim_k X_0 \( P_{n_k}(u)
\) $ exists, is a subgroup $ E $ of $ \R^\infty $, $ \g(E) = 1 $, and
the limit is the needed homomorphism $ X_1 : E \to G $.\qed

\medskip

The convergence $ X_0 \circ P_n + g \to X $, obtained in the proof
above, may be compared with other results [\JK, \Sa, \BI].
There, convergence almost sure is established for linear spaces;
here --- convergence in probability, for commutative groups.
Also, Lemma 5.3 may be compared with the study of ``quasi-additive
functionals'' in [\BI]. There, maps $ G \to \R $ are considered; here
--- maps $ \R^\infty \to G $.

\medskip

{\sc 5.4 Corollary.} Let $ G $ be (the additive group of) a separable 
F-space. Then, in Condition (c) of Lemma 5.3, the subgroup $ E $ can be
chosen to be a linear subspace, and the homomorphism $ X_1 $ --- a
linear map.

\medskip

{\sc Proof.}
The function $ X_0 : l_2 \to G $ is linear, since it is a continuous
homomorphism between F-spaces. Functions $ X_0 \circ P_{n_k} $
are linear, therefore their limit $ X_1 $ is linear, and its domain $
E $ is a linear subspace.\qed

\medskip

All reasonable definitions of a Gaussian measure on a separable Banach
space are evidently equivalent, which cannot be said about separable
F-spaces and Polish groups (see [\BI]).

A symmetric Gaussian measure in the sense of Fernique is a probability
measure $ \mu $ on a separable F-space $ F $ such that the product
measure $ \mu \otimes \mu $ on $ F \oplus F $ is invariant under the
following group of
transformations: $ T_\phi (x,y) = \( \, x \cos \phi - y \sin \phi, \,
x \sin \phi + y \cos \phi \, \) $ for $ x,y \in F $, $ \phi \in \R $.

A Gaussian measure in the sense of Bernstein is a probability measure
$ \mu $ on a commutative Polish group $ G $ such that the product
measure $ \mu \otimes \mu $ on $ G \times G $ turns into some product
measure $ \nu_1 \otimes \nu_2 $ on $ G \times G $ under the following
transformation: $ (x,y) \mapsto ( x-y, x+y ) $ for $ x,y \in G $.

If $ \mu $ is a symmetric Gaussian measure in the sense of Fernique,
then $ \mu $ (as well as any shift of $ \mu $) is also Gaussian in the
sense of Bernstein.

A Gaussian (convolution) semigroup is a family $
(\mu_t)_{t\in[0,\infty)} $ of probability measures $ \mu_t $ on a
Polish group $ G $, produced by some Brownian motion $ X(\cdot) $ in $
G $, in the sense that $ \mu_t $ is the distribution of $ X(t) $ for
each $ t $. (See [\BH].)

Let us define a {\it constructively Gaussian measure} as a probability
measure $ \mu $ on a commutative Polish group $ G $, that can be
represented as the distribution of some $ X \in L_0 (\g,G) $
satisfying (equivalent) conditions (a--c) of Lemma 5.3.

A constructively Gaussian measure is a Gaussian measure in the sense
of Bernstein. (Indeed, 5.3(c) transfers Bernstein property of $ \g $
into Bernstein property of $ \mu $.)

A constructively Gaussian measure on a separable F-space is a
shift of a symmetric Gaussian measure in the sense of Fernique. (Indeed,
5.3(c) and 5.4 transfer Fernique property of $ \g $ into Fernique
property of $ \mu $ shifted by $ (-g) $.)

If a measure on a separable F-space is contained in a Gaussian
semigroup, then it is constructively Gaussian (which is shown below by
means of Theorem 1.8). The following fact is thus obtained.

\medskip

{\sc 5.5 Corollary.}
If a measure $ \mu $ on a separable F-space is contained in a Gaussian
semigroup, then some shift of $ \mu $ is a symmetric Gaussian measure
in the sense of Fernique, and $ \mu $ is a Gaussian measure in the
sense of Bernstein.

\medskip

A Brownian motion in a separable F-space $ F $ determines a
measure $ \mu $
on the space $ C_0 \( [0,\infty), F \) $ of all continuous functions $
x : [0,\infty) \to F $ such that $ x(0)=0 $. The space $ C_0 \(
[0,\infty), F \) $, equipped with the topology uniform on finite
intervals, is also a separable F-space.

\medskip

{\sc Corollary 1.10} {\it(final formulation).\/}
For every Brownian motion in a separable F-space $ F $, the
corresponding measure $ \mu $ on the space $ C_0 \( [0,\infty), F \) $
is constructively Gaussian.

\medskip

The corresponding Gaussian semigroup $ (\mu_t) $ results from $ \mu $
by applying evaluation maps $ x(\cdot) \mapsto x(t) $. Each evaluation
map is a continuous linear map $ C_0 \( [0,\infty), F \) \to F $,
therefore it sends a constructively Gaussian measure into another
constructively Gaussian measure. Thus, 5.5 follows from 1.10.

\medskip

{\sc Proof of Corollary 1.10.} Let $ X $ be a Brownian motion in a
separable
F-space $ F $. By Theorem 1.8, there exists a Brownian motion $
(X,Y) $ in $ F \times \D $, where $ \D $ is the Hilbert space, such that
$ \F^X_t = \F^Y_t $ for all $ t $. Take a Borel function $ f : C_0 \(
[0,\infty), \D \) \to C_0 \( [0,\infty), F \) $ such that $ f(Y) = X $
almost sure (here $ X $ is treated as a random variable $ \O \to C_0
\( [0,\infty), F \) $).

Introduce two independent copies $ (X_1,Y_1) $ and $ (X_2,Y_2) $ of
the Brownian motion $ (X,Y) $. It is easy to see that the process $
(X_1+X_2, Y_1+Y_2) $ is also a Brownian motion. Moreover, the process
$ t \mapsto \( X_1(t/2)+X_2(t/2), Y_1(t/2)+Y_2(t/2) \) $ is another
copy of $ (X,Y) $, since $ \( X_1(t/2), X_2(t/2) \) $ is distributed
like $ \( X(t/2), X(t) - X(t/2) \) $, and the same for $ (X,Y) $
pairs. For notational convenience, introduce $ h : [0,\infty) \to
[0,\infty) $ by $ h(t) = t/2 $; we see that $ ( X_1 \circ h + X_2
\circ h, Y_1 \circ h + Y_2 \circ h ) $ is distributed like $ (X,Y) $.
Therefore, $ f (Y_1 \circ h + Y_2 \circ h) =  X_1 \circ h + X_2 \circ
h $ almost sure. However, $ f(Y_1) = X_1 $ and $ f(Y_2) = X_2 $, thus,
$ f (Y_1 \circ h + Y_2 \circ h) = f (Y_1) \circ h + f (Y_2) \circ h $.
Define $ g : C_0 \( [0,\infty), \D \) \to C_0 \( [0,\infty), F \) $ by
$ g(y) = f \( 2^{1/2} y \circ h ) \circ h^{-1} $, then $ g \( 2^{-1/2}
(Y_1+Y_2) \) = f(Y_1) + f(Y_2) $ almost sure. It means that $ f $
satisfies Condition (a) of Lemma 5.3. Though, the corresponding
measure $ \g $ is not the standard Gaussian measure on $ \R^\infty $,
rather it is some Gaussian measure on $ C_0 \( [0,\infty), \D \) $.
However, it is well-known that each Gaussian measure on a {\it locally
convex} F-space is linearly isomorphic $ \rm mod \, 0 $ to the
standard Gaussian measure on $ \R^\infty $.
Due to Lemma 5.3, $ f $ satisfies its Condition (c), therefore the
distribution of $ f(Y) = X $ is constructively Gaussian.\qed

\medskip

{\sc 5.6 Note.} Linearity of $ F $ was not used (multiplications by $
2^{\pm1/2} $ were made in $ \D $, not in $ F $). Thus, Corollary 1.10
(and its proof) remains true for all commutative Polish groups.

\bigskip
\centerline{\sc Appendix}
\nobreak\medskip

\def\odo#1#2{{\langle\!\langle#1,#2\rangle\!\rangle}}

{\sc Proof of Lemma 3.1.}
{\tolerance=10000
First, let $ H $ be a Hilbert space of finite or countable dimension,
$ X $ a Brownian motion in the unitary group $ \U (H) $, and $ (\T_t)
$ the semigroup on $ V $ defined by (3.3). We extend $ (\T_t) $ to the
complexification of $ V $, the space of all (not only Hermitian)
trace-class operators, by complex linearity; still, (3.3) holds.
Introduce a notation for
one-dimensional operators on $ H $: for any $ h_1, h_2 \in H $
$$  \odo{h_1}{h_2} : H \to H \, , \qquad
  \odo{h_1}{h_2} h = (h,h_2) h_1 \hbox{ for } h \in H \, .  $$
(Physicists denote it by $ |h_1\rangle\langle h_2| $.)
Note some general rules:
$$\eqalign{
  & \Tr \( \odo{h_1}{h_2} \) = (h_1,h_2) \> , \cr
  & \( \odo{h_1}{h_2} \)^* = \odo{h_2}{h_1} \> , \cr
  & A \odo{h_1}{h_2} = \odo{A h_1}{h_2} \, , \quad
    \odo{h_1}{h_2} A = \odo{h_1}{A^* h_2} \quad
    \hbox{for } A : H \to H \> , \cr
  & \odo{h_1}{h_2} \odo{h_3}{h_4} = (h_3,h_2) \odo{h_1}{h_4} \> , \cr
  & \Tr \( A \odo{h_1}{h_2} \) = (A h_1,h_2) \> . \cr
}$$
Using the Hilbert-Schmidt scalar product $ (B,C)_\HS = \Tr(BC^*) $ for
trace class
operators $ B,C : H \to H $, calculate a matrix element for $ \T_t $:
$ \( \T_t \odo{h_1}{h_2}, \odo{h_3}{h_4} )_\HS = \E \Tr \( X(t)
\odo{h_1}{h_2} X^*(t) \odo{h_4}{h_3} \) = \E \Tr \( \odo{X(t)
h_1}{X(t) h_2} \odo{h_4}{h_3} \) = \E \(h_4, X(t)h_2\) \(X(t)h_1,
h_3\) = \E \( X(t), \odo{h_3}{h_1} \)_\HS \( \odo{h_4}{h_2}, X(t)
\)_\HS $. In particular,
$$  \( \T_t \odo{h_1}{h_2}, \odo{h_1}{h_2} )_\HS = \E \( X(t), \rho_1
  \)_\HS \( \rho_2, X(t) \)_\HS \> ,  \leqno\rm(A.1)  $$
where $ \rho_k = \odo{h_k}{h_k} $.

}

From now on, $ H $ is assumed to be finite-dimensional. We have
diffusion processes $ A(t) $ in $ \D $ and $ X(t) = \exp (iA(t)) $ in
$ \U(H)
$; here $ \D $ is the linear space of all Hermitian operators on $ H $.
We may write $ A(t) = A_0 (t) + \l (t) \cdot 1_H $, $ \Tr A_0 (t) = 0
$, $ \l (t) \in \R $. For small $ t $,
$$\eqalign{
  & X(t) = 1 + i A(t) - \frt12 A^2(t) + o(t) \> ; \cr
  & \( X(t), \rho \)_\HS = 1 + i \( A(t), \rho \)_\HS - \frt12 \(
    A^2(t), \rho \)_\HS + o(t) \cr
}$$
for $ \rho = \odo h h $, $ \| h \| = 1 $. Therefore, for $ \rho_k =
\odo{h_k}{h_k} $, $ \| h_k \| = 1 $ ($ k=1,2 $),
$$\eqalign{
  & \E \( X(t), \rho_1 \)_\HS \( \rho_2, X(t) \)_\HS = \cr
  & = 1 + i \( \E A(t), \rho_1 \)_\HS - i \( \rho_2, \E A(t) \)_\HS +
    \cr
  & \quad + \E \( A(t), \rho_1 \)_\HS \( \rho_2, A(t) \)_\HS - \frt12
    \( \E A^2(t), \rho_1 \)_\HS - \frt12 \( \rho_2, \E A^2(t) \)_\HS
    + o(t)  \> . \cr
}$$
Separating real and imaginary terms (you see, $ A(t) $ and $ \rho_k $
are Hermitian) and using (A.1) we conclude that the following two
expressions are uniquely determined by the semigroup $ (\T_t) $:
$$\leqalignno{
  & \frac{d}{dt} \Big|_{t=0} \( \E A(t), \rho_1 - \rho_2 \)_\HS \> ;
    &\rm(A.2) \cr
  & \frac{d}{dt} \Big|_{t=0} \E \( A(t), \rho_1 \)_\HS \( A(t),
    \rho_2 \)_\HS - \frac12 \frac{d}{dt} \Big|_{t=0} \( \E A^2(t),
    \rho_1+\rho_2 \)_\HS \> . &\rm(A.3) \cr
}$$
Operators of the form $ \rho_1 - \rho_2 $ belong to the space $ \D_0 $
of all traceless Hermitian operators, and $ \D_0 $ is spanned by such
operators. Thus, (A.2) means that the following is uniquely determined
by $ (\T_t) $:
$$  \frac{d}{dt} \Big|_{t=0} \E A_0 (t) \> .  \leqno\rm(A.4)  $$
Consider four expressions of the form (A.3): first, exactly as (A.3),
that is, for the pair $ (\rho_1,\rho_2) $; second, the same for
another pair, $ (\rho_3,\rho_4) $; also, for $ (\rho_1,\rho_4) $ and $
(\rho_3,\rho_2) $, the last two with the minus sign. Summing the four,
we get
$$  \E \( A(t), \rho_1-\rho_3 \)_\HS \( A(t), \rho_2-\rho_4 \)_\HS \>
  .  \leqno\rm(A.5)  $$
It means that the whole bilinear form describing the spread (at $ t=0+
$) of $ A_0 (t) $ in $ \D_0 $ is determined by $ (\T_t) $. Together
with (A.4) it proves that the generator of the Brownian motion $ \(
\det X(t) \)^{-1/d} X(t) = \exp \( i A_0 (t) \) $ in $ \SU(d) $ is
uniquely determined by $ (\T_t) $, which completes the proof of
3.1(a).

Turn to the semigroup $ (x_t) $, $ x_t = \E X(t) $. We have $ x_t = 1
+ i \E A(t) - \frac12 \E A^2(t) + o(t) $; separating real and
imaginary terms we conclude that the following two expressions are
uniquely determined by $ (x_t) $:
$$\leqalignno{
  & \frac{d}{dt} \Big|_{t=0} \E A(t) \> , &\rm(A.6) \cr
  & \frac{d}{dt} \Big|_{t=0} \E A^2 (t) \> . &\rm(A.7) \cr
}$$
However, $ A(t) = A_0 (t) + \l(t) \cdot 1_H $, and $ \l(t) = (1/d) \Tr
A(t) $. Taking the trace of (A.6), we find $ \frac{d}{dt} \big|_{t=0}
\E \l(t) $. Further, $ \E A^2(t) = \E A_0^2(t) + 2 \E \l(t) A_0(t) +
\E \l^2(t) \cdot 1_H $, and 3.1(a) ensures that $ \E A_0^2(t) $ is
determined by $ (\T_t) $. Thus, the following is uniquely determined
by $ (x_t) $ and $ (\T_t) $:
$$  2 \frac{d}{dt} \Big|_{t=0} \E \l(t) A_0(t) + 1_H \cdot
  \frac{d}{dt} \Big|_{t=0} \E \l^2(t) \> .  \leqno\rm(A.8)  $$
The two terms can be separated by taking the trace. So, $ (x_t) $ and
$ (\T_t) $ determine
$$  \frac{d}{dt} \Big|_{t=0} \E \l(t) \, , \quad
  \frac{d}{dt} \Big|_{t=0} \E \l^2 (t) \, , \quad
  \frac{d}{dt} \Big|_{t=0} \E \l(t) A_0(t) \, ;  \leqno\rm(A.9)  $$
these are all the infinitesimal characteristics of $ A $ (or $ X $),
besides characteristics of $ A_0 $ given by 3.1(a).\qed

\bigskip
\centerline{\sc References}
\nobreak\medskip

\firstfalse

\bye